\documentclass[a4paper, reqno, 10pt]{amsart}
\usepackage{amssymb}
\usepackage{extarrows}
\usepackage{bm}
\usepackage[centering]{geometry}
\usepackage{enumitem}
\usepackage{longtable}
\usepackage{multirow}
\usepackage{array}
\usepackage{etex}
\usepackage{pictex}
\usepackage{tikz}
\usepackage{graphics}

\usepackage[all]{xy}

\usepackage{amsmath,amsthm}

\usepackage{pictex}

\DeclareMathOperator{\Ext}{Ext}

\newcommand{\ssize}[1]{\smallmatrix #1\endsmallmatrix}
\newcommand{\arr}[2]{\arrow <1.5mm> [0.25,0.75] from #1 to #2}

\numberwithin{equation}{section}
\theoremstyle{definition}
\newtheorem{defn}{Definition}[section]
\newtheorem{rem}[defn]{Remark}

\newtheorem{exm}[defn]{Example}

\theoremstyle{plain}
\newtheorem{cor}[defn]{Corollary}
\newtheorem{thm}[defn]{Theorem}
\newtheorem{lem}[defn]{Lemma}

\newtheorem{prop}[defn]{Proposition}
\newtheorem{fact}[defn]{Fact}

\newtheorem*{Freyd-Mitchell}{Freyd-Mitchell embedding Theorem}
\def\G{\text{\rm G}}

\def\op{\text{\rm op}}
\def\Coker{\operatorname{Coker}}

\def\gp{\operatorname{gp}}
\def\mod{\operatorname{mod}}

\def\Hom{\operatorname{Hom}}
\def\End{\operatorname{End}}
\def\Ext{\operatorname{Ext}}

\def\add{\operatorname{add}}
\def\Ker{\operatorname{Ker}}

\def\soc{\operatorname{soc}}
\def\Tr{\operatorname{Tr}}
\def\Im{\operatorname{Im}}

\def\bdim{\operatorname{\bold{dim}}}

\def\arr#1#2{\arrow <1.5mm> [0.25,0.75] from #1 to #2}

\def  \Cal{\mathcal }
\def\top{\operatorname{top}}
\def\soc{\operatorname{soc}}

\title[Claus Michael Ringel's main contributions to Gorenstein-projective modules]{Claus Michael Ringel's main contributions \\ to Gorenstein-projective modules}
\dedicatory{\bf Dedicated to Professor Ringel on the occasion of his eightieth birthday}
\author[Nan Gao,  Xue-Song Lu, Pu Zhang] {Nan Gao,  Xue-Song Lu, Pu Zhang$^*$}
\thanks{$^*$ Corresponding author}
\thanks{nangao@shu.edu.cn \ \ \ \ leocedar@sjtu.edu.cn \ \ \ \ pzhang$\symbol{64}$sjtu.edu.cn}
\thanks{Supported by National Natural Science Foundation of China $($Grant No. 12131015, 12271333$)$ and Natural Science Foundation of Shanghai $($Grant No. 23ZR1435100$)$.}
\begin{document}
\maketitle
\begin{abstract} In this article we try to recall Claus Michael Ringel's works on the Gorenstein-projective modules.
This will involve but not limited to his fundamental contributions, such as in,
the solution to the independence problem of totally reflexivity conditions;
the technique of $\mho$-quivers;
a fast algorithm to obtain the Gorenstein-projective modules over the Nakayama algebras;
the one to one correspondence between the indecomposable non-projective perfect differential modules of a quiver  and
the indecomposable representations of this quiver;
the description of the module category of the preprojective algebras of type $\mathbb A_n$ via submodule category;
semi-Gorenstein-projective modules, reflexive modules, Koszul modules, as well as the $\Omega$-growth of modules, over short local algebras;
and his negative answer to the question whether an algebra has to be self-injective in case all the simple modules are reflexive.

\vskip5pt


\end{abstract}

\vskip10pt

Since 2012 Claus Michael Ringel has paid attention to the Gorenstein-projective modules.
With his wide insights, deep thoughts, keen intuitions, and ingenious quiver skills, during 2012 - 2023,  he
has made a series of fundamental contributions.  The aim of this article is to try to
recall Professor Ringel's main works on the Gorenstein-projective modules.

\vskip5pt

This class of modules has been introduced  by Maurice Auslander [A] under the name {\it modules of Gorenstein dimension zero} in 1967, in his four lectures at the S\'eminaire Pierre Samuel.
Auslander discussed the class $G(A)$ of finitely generated modules satisfying the conditions (G1), (G2), (G3)
([A, Definition 3.2.2]), over a two-sided noetherian ring $A$.
In [AB, 3.8] M. Auslander and M. Bridger showed that a module $M\in G(A)$ if and only if $M$ and $\Tr M$ are semi-Gorenstein-projective.
E. E. Enochs and O. M. G. Jenda [EJ1] have defined {\it Gorenstein projective modules}
(not necessarily finitely generated) over any ring, via complete projective resolutions.
It is then known that for finitely generated modules over a two-sided noetherian ring, the class of Gorenstein projective modules
coincides with $G(A)$ ([Chr, Theorem 4.2.6]. See Theorem \ref{fundamentalG}).

\vskip5pt

If no otherwise stated, throughout $A$ is an artin algebra, and modules are finitely generated left $A$-modules.
Let $A$-$\mod$ be the category of $A$-modules
and $\add(A)$ the full subcategory of all projective modules. Denote by $\gp(A)$ the full subcategory of Gorenstein-projective modules.

\vskip5pt

The class $\gp(A)$ is closed under direct summands, extensions, and the kernel of epimorphisms ([AR2]).
The Gorenstein-projective modules have more stable property than projective modules ([AB], [SWSW]),
and they and reflexive modules play important roles in resolutions of singularities ([V]).
A striking feature is that $\gp(A)$ is one of the most important Frobenius categories ([Bel2], [Chen2]), and
the corresponding stable category $\underline{\gp}(A)$ is triangulated ([Hap1]).
If $A$ is Gorenstein then $\gp(A)$ has Auslander-Reiten sequences, and $\underline{\gp}(A)$ has Auslander-Reiten triangles, which is just the singularity category of $A$ ([Buch], [Hap2], [O]).
Gorenstein-projective modules are also called {\it maximal Cohen-Macaulay modules} ([Buch]), {\it totally reflexive modules} ([AM]),
and {\it Cohen-Macaulay modules} ([Bel2]).
They have wide applications in representation theory,
relative homological algebra, algebraic geometry (especially, in singularity theory), and so on
(see e.g. [Buch], [BGS], [AR2], [EJ2], [AM], [Hol],  [CFH], [IY], [CPST], [Bel3]).

\vskip5pt

It becomes a fundamental problem whether the conditions {(G1)}, {(G2)}, and {(G3)} in Theorem \ref{fundamentalG} are independent, after they were proposed
by Auslander as the first definition of the Gorenstein-projective modules.
L. L. Avramov and A. Martsinkowsky [AM] explicitly state this problem.
D. A. Jorgensen and L. M. \c Sega [JS] present the first examples of modules over a 8-dimensional commutative algebra satisfying (G1) and (G3) but not (G2),
and modules satisfying (G2) and (G3) but not (G1). However, modules satisfying (G1) and (G2) but not (G3)
are  more difficult to obtain, since such modules will induce other two kinds of modules.
C. M. Ringel proposed to look at a 6-dimensional non-commutative local algebra, and in [RZ4] it is proved that this algebra admits modules satisfying (G1) and (G2) but not (G3). Hence
Theorem \ref{G1G2notG3} (also Theorem \ref{independence}) completely answers the independence problem positively.
Thus,  bi-semi-Gorenstein-projective modules and weakly Gorenstein algebras are really of interest, and deserve to be studied independently. See Section 1.

\vskip5pt

C. M. Ringel realized the importance of $\mho$-sequence and $\mho$-quiver is introduced in [RZ4].
The operator $\mho$ is a kind of inverse of the syzygy $\Omega$, it coincides with
$\Tr\Omega\Tr$ which has been studied by Auslander and Reiten [AR4], and turns out to be a useful
tool in studying various modules.  See Section 2.
	
\vskip5pt

C. M. Ringel emphasizes roles of examples. However, it seems that at that time there were not so many examples of Gorenstein-projective modules.
So he looked at the Nakayama algebras $A$. With strong skills on quivers,
he could describe the structure of $\gp(A)$. By using resolution quiver,
he provides a fast algorithm to obtain the Gorenstein-projective modules, and to judge whether $A$ is Gorenstein or CM-free.
Usually $\gp(A)$ is not abelian, quite surprising, after smartly considering the Gorenstein core $\mathcal{C}$,
he even proves that $\mathcal{C}$ is equivalent to the module category of a self-injective Nakayama algebra.
Thus, $\mathcal{C}$ is an abelian Frobenius category with stable category  just $\underline{\gp}(A)$.
This becomes his first paper on the Gorenstein-projective modules. See Section 3.

\vskip5pt

For finite acyclic quiver $Q$ and $k$-algebra $A$, $\Lambda$-modules are just the representations of $Q$ over $A$, where $\Lambda = A\otimes_kkQ$.
The Gorenstein-projective $\Lambda$-modules can be described via the monomorphism category of the Gorenstein-projective $A$-modules.
However, it is hard to determine the indecomposability of a Gorenstein-projective $\Lambda$-module.
C. M. Ringel suggests to look at the representations of $Q$ over $A = k[x]/\langle x^2\rangle$.
In this case, $\Lambda$-modules are exactly the differential $kQ$-modules, and the Gorenstein-projective $\Lambda$-modules
are the perfect differential $kQ$-modules. Using this identification and the covering theory, very surprising, C. M. Ringel establishes a remarkable one to one correspondence
between indecomposable non-projective Gorenstein-projective $\Lambda$-modules and  indecomposable $kQ$-modules. See Theorems \ref{RZ3.1} and \ref{RZ3.2}.

\vskip5pt

Preprojective algebras have been introduced by I. M. Gelfand and V. A. Ponomarev [GP], and studied by many mathematicians,
for examples, by G. Lusztig [L] (Lagrangian variety) in quantum groups, and by C. Geiss, B. Leclerc and J. Schr\"oer [GLS] in cluster algebras.
On the other hand, submodule categories have been studied by G. Birkhoff [Bir]. They are developed as various monomorphism categories, and have attracted a lot of attention
(e.g. [RS1], [S1], [RS2], [RS3],  [Z1],  [Chen1], [KLM1], [KLM2], [S2], [ZX], [AHS], [Z2], [GKKP]).
C. M. Ringel has smartly observed that the module category of the preprojective algebra $\Pi_{n-1}$ is a factor category of the
submodule category $\mathcal{S}(n)$, via a dense, full and objective functor $\mathcal{S}(n) \longrightarrow \Pi_{n-1}$-mod. See Theorem \ref{RZ1.1}.
While $\mathcal{S}(n) = \gp(T_2(\Lambda_n))$,
where $\Lambda_n = k[x]/\langle x^n\rangle$. Since  $\Pi_{n-1}$-mod is also a factor category of $\underline {\mathcal{S}(n)}$,
this produces an abelian category from a triangulated category.

\vskip5pt

Short local algebra $A$ has a lot interest in commutative algebra (e.g. [Les], [AGP], [HSV], [CV], [AIS], [HJS]).
C. M. Ringel stresses its special role in representation theory of algebras.
Using the Hilbert-type $(e, a)$ and ${\bf dim} M$ of module $M$ with Loewy length at most $2$,
he could give a numerical necessary condition for the existence of non-projective reflexive $A$-modules  over non self-injective short local algebra $A$.
The least dimension of such an algebra is $6$ of Hilbert-type $(3, 2).$
Necessary conditions for the existences of minimal exact complexes of projective modules, of semi-Gorenstein-projective modules,
and of $\infty$-torsionfree modules, are also given. The Auslander-Reiten conjecture  which claims that if $M$ is non-projective
semi-Gorenstein-projective then $\Ext^i_A(M,M)\neq 0$ for some $i\ge 1$, is proved in a stronger form for short local algebras.
The $\Omega$-growth of modules and Koszul modules over $A$ are also studied.  See Section 6.

\vskip5pt

R. Marczinzik \cite{M3} has asked whether $A$ has to be self-injective in case all the simple modules are reflexive.
C. M. Ringel [R2] answered this question in the negative, by giving an $8$-dimensional wild algebra which is not self-injective, but the simple modules are reflexive.
See Example \ref{simplereflnotselfinj}.

\vskip5pt

Claus Michael Ringel belongs to the group of distinguished mathematicians who process exceptional qualities: they are full of creativity and energy of research; and if they enter a field, then the field will be greatly changed.

\vskip10pt

{\bf Contents}

1. \ Semi-Gorenstein-projective modules and Gorenstein-projective modules

\hskip15pt    1.1 \ Equivalent characterizations of the Gorenstein-projective modules

\hskip15pt    1.2 \ Bi-semi-Gorenstein-projective modules which are not torsionless

\hskip15pt    1.3 \ Independence problem

 \hskip15pt   1.4 \ 3-dimensional indecomposable $\Lambda(q)$-modules

 \hskip15pt   1.5 \ Left weakly  Gorenstein algebras

 \hskip15pt   1.6 \ Reduced bi-semi-Gorenstein-projective modules

2. \ The technique of $\mho$-quivers

3. \ Gorenstein-projective modules over the Nakayama algebras

4. \ Representations of quivers over the algebra of dual numbers

\hskip15pt   4.1 \ Perfect differential $kQ$-module

\hskip15pt   4.2 \ Indecomposable perfect differential $kQ$-modules

5. \ From submodule categories to preprojective algebras

\hskip15pt   5.1 \  Preprojective algebras of type $\mathbb A_n$

\hskip15pt   5.2 \  Objective triangle functors

6. \ Modules over short local algebras

\hskip15pt   6.1 \  Existence of non-projective reflexive modules

\hskip15pt   6.2 \ Short local algebras of dimension 6 with non-projective reflexive modules

\hskip15pt   6.3 \  Existence of exact complex of projective modules

\hskip15pt   6.4 \ Semi-Gorenstein-projective and $\infty$-torsionfree modules

\hskip15pt   6.5 \ The Auslander-Reiten conjecture

\hskip15pt   6.6 \ The $\Omega$-growth of modules

\hskip15pt   6.7 \ Koszul modules

\hskip15pt   6.8 \ Numerical descriptions of left Koszul algebras

7. \ Simple reflexive modules

\vspace {-8pt} \section{\bf Semi-Gorenstein-projective modules and Gorenstein-projective modules}
\subsection {Equivalent characterizations of the Gorenstein-projective modules}
An $A$-module $M$ is said to be {\it semi-Gorenstein-projective} provided that $\Ext^i_A(M, A) = 0$ for all $i\ge 1.$

\vskip5pt

Consider the $A$-dual $M^* = \Hom(M, A)$ of $M$. Then $M^*$ is a right $A$-module.
Let $\phi_M: M \longrightarrow M^{**}$ be the canonical $A$-map defined by $\phi_M(m)(f) = f(m)$ for $m\in M, \ f\in M^*$.
A module $M$ is {\it torsionless} if $\phi_M$ is injective, or, equivalently, if $M$ is a submodule of a projective module.
A module $M$ is {\it reflexive} if $\phi_M$ is bijective. Finitely generated projective modules are reflexive.

\vskip5pt

For a module $M$, let $PM$ be a projective cover of $M$, and
$\Omega M$ the kernel of the canonical map $PM \longrightarrow M$. The modules $\Omega^tM$ with $t\ge 0$
are called the {\it syzygy} modules of $M$.
A module $M$ is said to be {\it $\Omega$-periodic} if
 $\Omega^tM = M$  for some $t\ge 1$. Taking a minimal projective presentation $P(\Omega M) \stackrel f \longrightarrow  P M \longrightarrow M \longrightarrow 0$ of $M$,
then  one has an exact sequence $0\longrightarrow M^* \longrightarrow (P M)^*\stackrel{f^*}\longrightarrow (P(\Omega M))^* \longrightarrow \Coker f^*\longrightarrow 0$.
Put $\Tr M = \Coker f^*$, which is called  the {\it transpose of $M$}.

\vskip5pt

{\it A complete projective resolution} is a (double infinite) exact sequence
$$ P^\bullet: \qquad  \cdots \longrightarrow P^{-1}\longrightarrow P^{0}
  \overset {d^0} \longrightarrow P^{1}\longrightarrow \cdots$$
of projective left $A$-modules, such that $\Hom_A(P^\bullet, Q)$ is again exact for any projective module $Q$ (or, equivalently, $\Hom_A(P^\bullet, A)$ is exact).
A module $M$ is {\it Gorenstein-projective}, if there is a complete
projective resolution $P^\bullet$ such that $M$ is isomorphic to the image of $d^0.$ See [EJ1, EJ2].

\vskip5pt

Projective modules are Gorenstein-projective, and Gorenstein-projective modules are semi-Gorenstein-projective.
Denote by $\gp(A)$ the class of all left Gorenstein-projective modules
and by ${}^\perp A$ the class of all left semi-Gorenstein-projective modules. Then
$\gp(A) \subseteq {}^\perp A.$ The following result is fundamental and well-known (see e.g. [Chr, Theorem 4.2.6]).

\begin{thm} \label{fundamentalG} \ Let $A$ be an artin algebra, and $M\in A${\rm-mod}. Then the following are equivalent.

\vskip5pt

$(1)$ \  $M$  satisfies the conditions $(\G1), \ (\G2)$ and $(\G3):$

\hskip18pt ${\rm (G1)}$ \ $M$ is semi-Gorenstein-projective.

\hskip18pt ${\rm (G2)}$ \ $M^* = \Hom(M, \ _AA)$ is semi-Gorenstein-projective, as a right $A$-module.

\hskip18pt ${\rm (G3)}$ \ $M$ is reflexive.

\vskip5pt

$(2)$ \  $M$ is Gorenstein-projective.

\vskip5pt

$(3)$ \  $M$ and $\Tr M$ are semi-Gorenstein-projective.
\end{thm}

\begin{rem} \ (1) \ The equivalence of $(1)$ and $(3)$ is shown  in [AB, Proposition 3.8];
and the equivalence of $(1)$ and $(2)$ can be found in  [Chr, Theorem 4.2.6] where L. W. Christensen attributes the proof to L. L. Avramov, R.-O. Buchweitz, A. Martsinkovsky, and I. Reiten.
Note that $\Tr M$ is semi-Gorenstein-projective if and only if $M$ satisfies (G2) and (G3).

\vskip5pt

(2) \ Theorem \ref{fundamentalG} holds also for finitely generated modules over a two-sided noetherian ring. In this case,
although a minimal projective presentation of $M$ does not exist in general, one can replace $\Tr M$ by $\Coker f^*$, where
$P_1 \stackrel f \longrightarrow  P_0 \longrightarrow M \longrightarrow 0$ is a projective presentation of $M$.

\vskip5pt

(3) \ If $M$ is not finitely generated, then Theorem \ref{fundamentalG} does not hold.
\end{rem}

\subsection{Bi-semi-Gorenstein-projective modules which are not torsionless}
It is not easy to get semi-Gorenstein-projective modules which are not Gorenstein-projective.
The first such an example was constructed by Jorgensen and \c{S}ega [JS] in 2006,
for a commutative algebra of dimension 8. In 2017 Marczinzik [M2] presented such modules over some non-commutative algebras.
It is more difficult to get bi-semi-Gorenstein-projective modules $M$ which are not Gorenstein-projective,
here by {\it a bi-semi-Gorenstein-projective module $M$} we mean that both $M$ and $M^*$ are semi-Gorenstein-projective.
To find bi-semi-Gorenstein-projective modules which are not Gorenstein-projective,
after many trials, C. M. Ringel proposed the following algebra.

\vskip5pt

Let $k$ be a field and $q$ a non-zero element in $k$.  Consider the $k$-algebra $\Lambda = \Lambda(q)$ of dimension $6$, which is generated by
$x, y, z$, with the
relations:
$$ x^2,\ y^2,\ z^2,\ yz,\ xy+qyx,\ xz-zx,\ zy-zx.$$
For $\alpha\in k$, let $M(\alpha)$ be the $\Lambda$-module with basis $v, v', v''$, such that $$xv = \alpha v', \ yv = v', zv = v'', xv' = xv'' = yv' = yv'' = zv' = zv'' =0.$$

\begin{thm} \label{G1G2notG3} \ {\rm([RZ4, Theorem 1.6])} \ If the multiplicative order $o(q)$ is infinite,
then $M(q)$ is bi-semi-Gorenstein-projective with $M(q)^*\cong (x-y)\Lambda$, but not torsionless,
thus $M(q)$ is not Gorenstein-projective.
	
\vskip5pt

If $\alpha\notin \langle q\rangle$, the cyclic group generated by $q$, then $\Lambda(q)$-modules
$M(\alpha)$ are Gorenstein-projective. In particular, $M(0)$ is Gorenstein-projective of $\Omega$-period $1$.
Thus, if $k$ is infinite, then there are infinitely many isomorphism classes of \ $3$-dimensional
Gorenstein-projective $\Lambda(q)$-modules.
\end{thm}

\subsection{Independence problem}
It was asked by L. L. Avramov and A. Martsinkowsky [AM] whether the
conditions {(G1)}, {(G2)}, {(G3)} in Theorem \ref{fundamentalG} are independent.
In [JS], Jorgensen and \c Sega present the first examples of modules satisfying (G1) and (G3), but not (G2),
and also modules satisfying (G2) and (G3), but not (G1). Combining this
with  Theorem \ref{G1G2notG3} one has
	
\vskip5pt

\begin{thm} \label{independence} \ {\rm([RZ4, Theorem 1.7])} \ For artin algebras,
the conditions $(\G1),\ (\G2)$ and $(\G3)$ are independent.

\vskip5pt

Moreover, if a module $M$ is semi-Gorenstein-projective and not Gorenstein-projective,
then $\Omega^2 M$ satisfies {\rm(G1)}  and {\rm(G3)}, but not {\rm(G2)}.

\vskip5pt

If a module $M$ satisfies {\rm(G1)} and {\rm(G3)}, but not {\rm(G2)}, then $M^*$ is a right module satisfying {\rm(G2)} and {\rm(G3)}, but not {\rm(G1)}.
\end{thm}
	
\begin{rem} \label{CW'semail} \ Theorem \ref{independence} first completely solves the independence problem, as claimed by L. W. Christensen and D. J. Wu, in the email to C. M. Ringel and P. Zhang on 2018-9-10:

\vskip5pt

``We just finished reading your paper Gorenstein-projective and semi-Gorenstein-projective modules."
``As best we can tell there is no example in [JS] of a module that satisfies (G1) and (G2) but not (G3), so your paper is probably the first one that ``truly"
demonstrates the independence of the three conditions."
\end{rem}

\begin{rem} \ Although R. Marczinzik [M2] gave semi-Gorenstein-projective modules which are not Gorenstein-projective, his result does not touch the independence problem.

\begin{thm} \ {\rm ([M2])} \ Let $A$ be a symmetric algebra, $M = A\oplus X$, where $X$ is an $A$-module without non-zero projective direct summand, and $B = {\rm End}_A(M)^{\rm op}$.
Suppose that $N$ is an $A$-module and $l$ is a positive integer such that $\Ext_A^l(X,N)\ne 0$ and $\Ext_A^i(X,N)=0$ for all $i\geq l+1$.
Let $R=\Hom_A(M,\Omega^{-l}N)\in B\mbox{-}\mod$. Then

$(1)$ \ $\Hom_B(D(B_B), \ _BB) \cong B$ as $B$-$B$-bimodules, where $D = \Hom_k(-, k)$.

$(2)$ \ ${\rm dom. dim} _BR=\infty$, which implies that $R$ is a semi-Gorenstein-injective $B$-module, i.e., $\Ext^i_B(D(B_B), R) = 0$ for all $i\ge 1.$

$(3)$ \ ${\rm codom. dim} _BD(B_B) \ge 1$. In particular, each left injective $B$-module is a quotient of projective-injective $B$-module.

$(4)$ \ ${\rm codom. dim} _BR=0$, which implies that $R$ is not a quotient of injective $B$-module. In particular, $D(R)$ is not torsionless. \end{thm}

Thus $D(R)$ satisfies ${\rm (G1)}$ but not ${\rm (G3)}$.  In [M2] the existences of such $A$, $l$,  $X$ and $N$,
can be guaranteed.  However, it is  not clear whether $D(R)$ satisfies ${\rm (G2)}$.
\end{rem}

\subsection{3-dimensional indecomposable $\Lambda(q)$-modules}

Let $k$ be a field and $q\in k-\{0\}$. Recall the algebra
$$\Lambda=\Lambda(q) = k\langle x, y, z\rangle / \langle
 x^2, \ y^2, \ z^2, \ yz, \ xy+qyx, \ xz-zx,\ zy-zx\rangle.$$
It is a short local algebra, here {\it short} means $J^3 = 0$, where $J$ is the Jacobson radical of $\Lambda$.
For $(a, b, c)\in k^3 - \{(0, 0, 0)\}$, consider $3$-dimensional left
 $\Lambda$-modules $$M(a,b,c) = {}_\Lambda\Lambda/(\Lambda(ax+by+cz) + \Lambda yx + \Lambda zx).$$
Then $M(a,b,c)\cong M(a', b' , c')$ as $\Lambda$-module if and only if $(a, b, c) = \lambda (a', b', c')$ for some $\lambda\in k^* = k-\{0\}$.
For example, for $\alpha\in k$, the $\Lambda$-module $M(\alpha)$ in Subsection 1.2 is isomorphic to
$$M(1, -\alpha, 0) = {}_\Lambda\Lambda/[\Lambda(x-\alpha y) + {\rm soc} \Lambda] = k\bar 1 \oplus \bar y \oplus \bar z.$$

Let $M$ be an indecomposable $3$-dimensional left $\Lambda$-module. Then $M$ is annihilated by  $J^2$,
and hence $M$ is a left $(\Lambda/J^2)$-module. Since $\Lambda/J^2$ is a commutative
 $k$-algebra, $D(M) = \Hom(M,k)$ is a left $(\Lambda/J^2)$-module, thus
 a left $\Lambda$-module.

\begin{prop} \label{3dimmod} \ {\rm([RZ5, A.1])} \ Let $M$ be an indecomposable $3$-dimensional left
 $\Lambda$-module. Then $M$ is isomorphic to one of the
 following pairwise non-isomorphic $\Lambda$-modules$:$
$$M(1, b, c) \ (b, c\in k^*), \ \ \ M(1, b, 0) \ (b\in k^*), \ \ \ M(1, 0, c) \ (c\in k^*),$$
$$M(0, 1, c) \ (c\in k^*), \ \ \ M(0, 0, 1), \ \ \ M(0, 1, 0), \ \ \ M(1, 0, 0),$$
$$D(M(1, b, c)) \ (b, c\in k^*), \ \ \ D(M(1, b, 0)) \ (b\in k^*), \ \ \ D(M(1, 0, c)) \ (c\in k^*),$$
$$D(M(0, 1, c)) \ (c\in k^*), \ \ \ D(M(0, 0, 1)), \ \ \ D(M(0, 1, 0)), \ \ \ D(M(1, 0, 0)).$$
\end{prop}

 \vskip5pt

\begin{thm} \label{detailedleft} \ {\rm([RZ5, Theorem 1.5])} \ Let $M$ be an indecomposable  left $\Lambda$-module of dimension $\le 3$. Then

\vskip5pt

$(1)$ \ $M$ is  torsionless if and only if $M$ is simple or $M\cong \Lambda(x-y)$, or $M\cong \Lambda z$, or $M\cong M(1,b,c)$ with $b\neq -q$, or $M\cong M(0,1,0)$, or $M\cong M(0,0,1).$

\vskip5pt

$(2)$ \ $M$ is  reflexive if and only if   $M\cong M(1, b, c)$, where $b\neq -q$ and $b\neq -q^2$.

\vskip5pt

$(3)$ \ $\Ext^1_\Lambda(M, \Lambda) = 0$  if and only if   \ $M\cong M(1, b, c)$, where $b\neq -1$.

\vskip5pt

$(4)$ \ $M$ is  semi-Gorenstein-projective  if and only if   \ $M\cong M(1, b, c)$, where $b\neq -q^i$ for $i\le 0.$

\vskip5pt

$(5)$ \ $M$ is semi-Gorenstein-projective and $M$ is not torsionless if and only if  $o(q) = \infty$ and $M\cong M(1, -q, c)$.

\vskip5pt

$(6)$ \ $\Tr M$ is semi-Gorenstein-projective if and only if  $M\cong M(1, b, c)$, where $b\neq -q^i$ for $i \ge 1.$

\vskip5pt

$(7)$ \ $\Tr M$ is semi-Gorenstein-projective and $\Ext^1_\Lambda(M, \Lambda)\ne 0$ if and only if  $o(q) = \infty$ and $M\cong M(1, -1, c)$.

\vskip5pt

$(8)$ \ $M$ is  Gorenstein-projective   if and only if   $M\cong M(1, b, c)$,  where $b\neq -q^i$ for $i\in \Bbb Z$.\end{thm}

\vskip5pt

Also, For $(a, b, c)\in k^3 - \{(0, 0, 0)\}$, consider $3$-dimensional right
$\Lambda$-modules
$$M'(a,b,c): = \Lambda_\Lambda/((ax+by+cz)\Lambda + yx \Lambda + zx \Lambda).$$

\begin{thm} \label{detailedright} \ {\rm([RZ5, Theorem 8.6])} \ Let $N$ be an indecomposable  right $\Lambda$-module of dimension $\le 3$. Then

\vskip5pt

$(1)$ \ $N$ is  torsionless if and only if  torsionless iff  $N$ is simple or $N\cong y\Lambda$, or $N\cong z\Lambda$, or $N\cong M'(1, b, c)$ with $b\neq -1$, or $N\cong M'(1,-1,0)$, or $N\cong M'(0,0,1).$

\vskip5pt

$(2)$ \ $N$ is  reflexive if and only if   $N\cong M'(1,b,c)$,  where $b\neq -q^{-1}$ and $b\neq - 1$.

\vskip5pt

$(3)$ \ $\Ext^1_\Lambda(N, \Lambda_\Lambda) = 0$ if and only if  $N\cong M'(1,b,c)$, where $b\neq -q.$

\vskip5pt

$(4)$ \ $N$ is semi-Gorenstein-projective  if and only if   $N\cong M'(1, b, c)$ with $b\neq -q^i$
for $i\ge 0$, or $N\cong M'(1, -1, c)$ with $c\ne 0.$

\vskip5pt

$(5)$ \ $N$ is semi-Gorenstein-projective and $N$ is not torsionless if and only if  $o(q) = \infty$ and $N\cong M'(1,-1,c)$ with $c\neq 0.$

\vskip5pt

$(6)$ \ $\Tr N$ is semi-Gorenstein-projective if and only if  $N\cong M'(1,b,c)$, where $b\neq -q^i$ for $i \le 0.$

\vskip5pt

$(7)$ \ $\Tr N$ is semi-Gorenstein-projective and $\Ext^1_\Lambda(N, \Lambda_\Lambda)\ne 0$ if and only if $o(q) = \infty$ and
   $N\cong M'(1,-q, c)$.

\vskip5pt

$(8)$ \ $N$ is  Gorenstein-projective if and only if  $N\cong M'(1, b, c)$, where $b\neq -q^i$ for $i\in \Bbb Z$. \end{thm}

One has key conversion formulas ([RZ5, Propositions 9.4, 9.6]):
$$M(1, b, c)^*\cong M'(1, q^{-2}b, (1+q^{-2}b)(1+q^{-1}b)c);$$ and \begin{align*}&M'(1, b, c)^* \cong M(1, q^2b, \frac {c}{(1+qb)(1+b)}),  \ \ \mbox{if} \  b\notin \{-1, -q^{-1}\};
\\ & M'(1, -1, c)^*\cong \Lambda z\oplus kyx, \ \ \mbox{if} \ c\ne 0;
\\ & M'(1, -1, 0)^*\cong \Lambda(x-qy) + \Lambda z;
\\ & M'(1, -q^{-1}, c)^*\cong \Lambda/(\Lambda z\oplus kyx) = M(0, 0, 1), \ \ \mbox{if} \ c\ne 0, \ q\ne 1;
\\ & M'(1, -q^{-1}, 0)^*\cong M(1, q, 0)\cong \Lambda(x-y)\Lambda, \ \ \mbox{if} \ q\ne 1.\end{align*}

With all these information one has

\begin{cor} \label{detailedconclusion} \ Let $o(q) = \infty$ and $M$ be an indecomposable $\Lambda$-module of dimension $\le 3$. Then

\vskip5pt

$(1)$ \ $M$ satisfies \ {\rm (G1)} \ and \ {\rm (G2)} \ but \ $M$ \ is not torsionless if and only if \  $M\cong M(1, -q, c)$ \ with \ $c\in k;$
if and only if $M$ satisfies \ {\rm (G1)} \ and \ {\rm (G2)} \ but \ $\varphi_M$ \ is not surjective.

\vskip5pt

For example, $M(q) \cong M(1, -q, 0)$ satisfies \ {\rm (G1)} \ and \ {\rm (G2)} \ but \ $M(q)$ \ is not torsionless.

\vskip5pt

$(2)$ \ $M$ satisfies \ {\rm (G1)} \ and \ {\rm (G3)} \ but not \ {\rm (G2)} if and only if \ $M\cong M(1, -q^i, c)$ with \ $i \ge 3, \ c\in k$.

\vskip5pt

$(3)$ \ $M$ satisfies \ {\rm (G2)} \ and \ {\rm (G3)} \ but not \ {\rm (G1)} if and only if \ $M\cong M(1,-q^i, c)$ with \ $i\le 0, \ c\in k$.

\vskip5pt

$(4)$ \   $M$ is Gorenstein-projective if and only if  $M\cong M(1, b, c)$ \ with \ $b\neq -q^i$ for all \ $i\in \Bbb Z$, \ $c\in k$.

\vskip5pt

$(5)$ \   $M$ is Gorenstein-projective if and only if  $M$ satisfies \ {\rm (G1)} \ and \ {\rm (G2)}, and \ $M$ \ is torsionless$;$ and
if and only if  $M$ satisfies \ {\rm (G1)} \ and \ {\rm (G2)}, and \ $\varphi_M$ \ is surjective. \end{cor}

\begin{rem} \label{problem2} $(1)$ \ By Corollary \ref{detailedconclusion}, the $\Lambda$-modules $M$ of dimension $\le 3$ satisfying \ {\rm (G1)} \ and \ {\rm (G3)} \ but not \ {\rm (G2)},
and the $\Lambda$-modules $M$ of dimension $\le 3$ satisfying \ {\rm (G2)} \ and \ {\rm (G3)} \ but not \ {\rm (G1)},  have two independent parameters; but
the $\Lambda$-modules $M$ of dimension $\le 3$ satisfying \ {\rm (G1)} \ and \ {\rm (G2)} \ but not \ {\rm (G3)}  have only one parameter.
Although this is the situation for the algebra $\Lambda$, it explains  in some sense why the
modules $M$ satisfying \ {\rm (G1)} \ and \ {\rm (G2)} \ but not \ {\rm (G3)} are more difficult to find out than the other two classes of modules. See also
Theorem \ref{independence}.

\vskip5pt

$(2)$ \ By Corollary \ref{detailedconclusion}, $\Lambda$-module $M$ of dimension $\le 3$ satisfying  {\rm (G1)}  and  {\rm (G2)}
is torsionless if and only if  $\varphi_M$  is surjective. In fact,
all the known bi-semi-Gorenstein-projective modules $M$ such that $\phi_M$ is an epimorphism or a
monomorphism, are Gorenstein-projective ([RZ6, 3.1]).
\end{rem}

\subsection {\bf Left weakly  Gorenstein algebras}

An artin algebra $A$ is {\it left weakly  Gorenstein}
if ${}^\perp A = \gp(A)$, i.e.,
any left semi-Gorenstein-projective module is Gorenstein-projective.

\vskip5pt

An artin algebra $A$ is {\it Gorenstein }
if the injective dimensions of $_AA$ and $A_A$ are finite. In this case the two dimensions are the same.
A Gorenstein algebra is left and right weakly  Gorenstein; but the converse is not true.
For example, the algebra given by the quiver $\xymatrix{\ar@(ur, ul)^\beta 2
\ar@{->}@<2pt> [r]^-{\alpha} &*-<1pt> {\begin{matrix} \ 1
  \end{matrix}}}$ and the relations $\beta^2, \ \alpha\beta$, is left and right weakly  Gorenstein, but not Gorenstein.

\vskip5pt

Left weakly Gorenstein algebras have been called ``nearly Gorenstein algebras" in [M2],
and  ``algebras with condition (GC)'' in [CH].
C. M. Ringel and his collaborator have apologized for using a new name, but they explained: ``the algebras with $\gp(A) = {}^\perp A$ seem
to be quite far away from being Gorenstein" ([RZ4, 1.9]).

\begin{thm} \label{weakG1} \ {\rm([RZ4, 1.2])} \ Let $A$ be an artin algebra. Then the  following are equivalent$:$

\vskip5pt

{\rm (1)} \ $A$ is left weakly  Gorenstein.

\vskip5pt

{\rm (2)} \ {\rm ([HH, Theorem 4.2])} \ Any semi-Gorenstein-projective module is torsionless.

\vskip5pt

{\rm (3)} \ Any semi-Gorenstein-projective module is reflexive.

\vskip5pt

{\rm (4)} \ For any semi-Gorenstein-projective module $M$, the map $\phi_M$ is surjective.

\vskip5pt

{\rm (5)} \ For any semi-Gorenstein-projective module $M$,  $M^*$ is semi-Gorenstein projective.

\vskip5pt

{\rm (6)} \ Any semi-Gorenstein-projective module $M$ satisfies $\Ext^1(M^*,A_A) = 0$.

\vskip5pt

{\rm (7)} \ Any semi-Gorenstein-projective module $M$ satisfies $\Ext^1(\Tr M,A_A) = 0.$\end{thm}
	
\vskip5pt

If $\Cal C$ is an extension-closed full subcategory of $A\mbox{-}\mod$, then
the embedding of $\Cal C$ into $A\mbox{-}\mod$ gives an exact structure on $\Cal C$.
An exact category $\Cal F$ is {\it a Frobenius category} if it has enough
projective and enough injective objects and the
projective objects of $\Cal F$ coincide with the injective objects of $\Cal F$ ([Kel1]).
Denote by $\Cal P(\Cal F)$ the full subcategory of the
projective objects in $\Cal F$.

\begin{thm} \label{weakG2} \ {\rm([RZ4, 1.3, 1.4])} \ Let $A$ be an artin algebra.

\vskip5pt

$(1)$ \  If the number of isomorphism classes of indecomposable modules which are both
semi-Gorenstein-projective and torsionless is finite, then $A$ is left weakly  Gorenstein
and any indecomposable non-projective semi-Gorenstein-projective module is $\Omega$-periodic.
	
\vskip5pt

$(2)$ \ Let $\Cal F$ be an  extension-closed full subcategory of $A\mbox{-}\mod$ such that
$\Cal F$ is a Frobenius category with respect to its canonical exact structure.
If $\Cal P(\Cal F) \subseteq {}^\perp A \subseteq \Cal F$, then
$A$ is left weakly Gorenstein and
$\Cal F = \gp(A).$

\vskip5pt

$(3)$ \ $A$ is left weakly  Gorenstein if and only if
${}^\perp\! A$ with its canonical exact structure is a Frobenius subcategory. \end{thm}

Y. Yoshino [Y] has shown
that for certain commutative rings $R$ the
finiteness of the number of isomorphism classes of indecomposable
semi-Gorenstein-projective $R$-modules implies that $R$ is left weakly  Gorenstein. This was generalized
to artin algebras by A. Beligiannis [Bel3, 5.11].
By R. Marczinzik [M1], all torsionless-finite artin algebras (i.e., there are only
finitely many isomorphism classes of torsionless indecomposable modules) are left weakly  Gorenstein.
Some classes of torsionless-finite algebras are given in [RZ4, 3.6].
Also, $\gp(A)$ is the largest resolving Frobenius subcategory of $\mod A$ ([Bel1, 2.11], [Bel2, p.145], and [Bel3, p.1989]; also [ZX, 5.1]).

\subsection{\bf Reduced bi-semi-Gorenstein-projective modules}

A bi-semi-Gorenstein-projective module is {\it reduced}
if it has no non-zero projective direct summands.
Such a  module is not necessarily Gorenstein-projective. Recall that  a complex of projective modules is {\it minimal} provided the image of each differential
is contained in the corresponding radical.

\begin{thm} \label{bisgp} {\rm ([RZ6, Main Theorem])} \  The isomorphism classes of reduced bi-semi-Gorenstein-projective $M$ correspond bijectively to the
isomorphism classes of minimal
complexes $P_\bullet$ of projective modules with ${\rm H}_i(P_\bullet) = 0$
for $i\neq 0,-1,$ so that the $A$-dual complex $P_\bullet^*$ is acyclic.

\vskip5pt

Explicitly, let $M$ be a reduced bi-semi-Gorenstein-projective module. Taking minimal projective resolutions
$$\cdots \longrightarrow P_2 \stackrel {f_2} \longrightarrow P_1  \stackrel {f_1} \longrightarrow P_0  \stackrel {e} \longrightarrow M \longrightarrow 0 \qquad \text{and} \qquad
 0 \longleftarrow  M^* \stackrel {q}\longleftarrow Q_0 \stackrel {d_1}\longleftarrow Q_1 \stackrel {d_2}\longleftrightarrow  Q_2 \longleftarrow \cdots$$
and for $i < 0$, let $P_i = Q_{-i-1}^*$ and $f_i = d_{-i}^*: P_{-i} \longrightarrow P_{-i-1}$,
and let $f_0 = q^*\phi_Me: P_0 \longrightarrow P_{-1}$. In this way, one obtains a minimal complex
$P_\bullet$ of projective modules
$$
\cdots \longrightarrow P_1 \stackrel {f_1} \longrightarrow P_0 \stackrel {f_0} \longrightarrow P_{-1} \stackrel {f_{-1}} \longrightarrow P_{-2} \longrightarrow \cdots$$
with ${\rm H}_i(P_\bullet) = 0$ for $i\neq 0,-1,$ such that the $A$-dual complex $P_\bullet^*$ is acyclic.
By construction, $M = \Coker f_1$, \  $M^* = \Coker d_1 = \Coker f_{-1}^*$, and
$$
 H_{0}(P_\bullet) = \Ker \phi_M \quad\text{\it and}\quad
 H_{-1}(P_\bullet) = \Coker \phi_M.$$

Conversely, let \ $P_\bullet = (P_i,f_i: P_i\longrightarrow P_{i-1})_i$
be a minimal complex of projective modules
with ${\rm H}_i(P_\bullet) = 0$
for $i\neq 0,-1$, so that $P_\bullet^*$ is acyclic.
Then $M = \Coker f_1$ is a reduced bi-semi-Gorenstein-projective module with $M^* = \Coker f_{-1}^*$ and $M^{**} = \Ker f_{-1}$.\end{thm}

\vskip5pt

Also, the transpose $\Tr$ provides a bijection between reduced
bi-semi-Gorenstein-projective modules and
reduced $\infty$-torsionfree right modules $Z$
with $\Omega^2Z$ semi-Gorenstein-projective ([RZ6]).

\section {\bf The technique of $\mho$-quivers}

C. M. Ringel realized the importance of $\mho$-sequence and $\mho$-quiver is introduced in [RZ4]. This turn out to be a useful
tool in studying various modules.
Denote by $\mho M$ the cokernel of a minimal left $\add(A)$-approximation of $M$.
The symbol $\mho$, pronounced
``agemo'', is a reminder that it is a kind of inverse of $\Omega$.
It turns out that the operator $\mho$ coincides with $\Tr\Omega\Tr$ ([AR4]).

\vskip5pt

Let $M$ be an indecomposable non-projective $A$-module, and $\omega: M \longrightarrow P$ a minimal left $\add(A)$-approximation with cokernel map
$\pi: P \to \mho M$. Then the image of $\omega$
is contained in the radical of $P$, thus $\pi$ is a projective cover.
If $M$ is in addition torsionless,  then $\omega$ is injective and $\mho M$ is indecomposable
and not projective, and $\Omega \mho M \simeq M.$

\vskip5pt

An exact sequence $0 \longrightarrow M \stackrel \omega \longrightarrow P \longrightarrow Z \longrightarrow 0$ in $A$-mod is an
{\it $\mho$-sequence} provided that both $M$ and $Z$ are indecomposable and not projective, and $\omega$ is a left $\add(A)$-approximation; or equivalently,
$M$ is indecomposable and not projective, and $\omega$ is a minimal left $\add(A)$-approximation.
Such a $\mho$-sequence exists if and only if $M$ is indecomposable torsionless and not projective.

\vskip5pt

The {\it $\mho$-quiver} of $A$ has as vertices the isomorphism classes $[X]$ of the
indecomposable non-projective modules $X$ and there is an arrow
$$\xymatrix {[X] & [\mho X]\ar@{-->}[l]}$$
if $X$ is torsionless. (The orientation of the arrow corresponds to the usual one for the  $\Ext$-quiver.)
Thus at any vertex of the $\mho$-quiver, at most one arrow starts and at most one arrow
ends.
	
\vskip5pt

Any $\mho$-component is a linearly oriented quiver
$\Bbb A_n$ with $n\ge 1$ vertices, or an oriented cycle $\widetilde{\Bbb A}_n$ with
$n\!+\!1\ge 1$ vertices, or of the form $-\Bbb N,$ or $\Bbb N,$  or $\Bbb Z$,
here are the quivers $-\Bbb N$ and $\Bbb N$:
$$
{\beginpicture
\setcoordinatesystem units <1cm,1cm>
\put{\beginpicture
\multiput{$\circ$} at 1 0  2 0  3  0 /
\setdashes <1mm>
\arr{0.8 0}{0.2 0}
\arr{1.8 0}{1.2 0}
\arr{2.8 0}{2.2 0}
\put{$\cdots$} at -.3 0
\put{$-\Bbb N$} at 1.5 -.5
\endpicture} at 0 0
%%%%%%%%%%%%%%%%%%%%%%
\put{\beginpicture
\multiput{$\circ$} at 0 0  1 0  2 0 /
\setdashes <1mm>
\arr{0.8 0}{0.2 0}
\arr{1.8 0}{1.2 0}
\arr{2.8 0}{2.2 0}
\put{$\cdots$} at 3.3 0
\put{$\Bbb N$} at 1.5 -.5
\endpicture} at 5 0
\endpicture}
$$
and all cases can arise as $\mho$-components. See [RZ4, 1.5, 7.7].
An indecomposable non-projective module $M$ will be said to be of {\it $\mho$-type}
$\Delta$ where
$\Delta \in \{\Bbb A_n,\ \widetilde{\Bbb A}_n,\, -\Bbb N,\ \Bbb N,\ \Bbb Z\}$
in case the $\mho$-component containing $[M]$ is of the form $\Delta$.
	
\vskip5pt

Many information on modules can be read out from its position in the $\mho$-quiver.
A module $M$ is  {\it $t$-torsionfree} if $\Ext^i(\Tr M,A_A) = 0$
for $1\le i \le t$ (and {\it $\infty$-torsionfree,} if $\Ext^i(\Tr M,A_A) = 0$
for all $i \ge 1$, i.e., $\Tr M$ is semi-Gorenstein-projective)  (see [AB]).

\begin{thm} \label{quiver1} \ {\rm([RZ4, 1.5])} \ Let $M$ be an indecomposable non-projective module.
	
\vskip5pt

{\rm(0)} \ $[M]$ is an isolated vertex \ if and only if \ $\Ext^1(M,A) \neq 0$ and $M$ is not torsionless.
	
\vskip5pt

{\rm(1)} \ $[M]$ is the start of a path of length $t\ge 1$ \ if and only if \
  $\Ext^i(M,A) = 0$ for $1\le i \le t.$

  In particular, $[M]$ is the start of an arrow \ if and only if\  $\Ext^1(M,A) = 0.$

\vskip5pt

{\rm(1$'$)} \ $[M]$ is the start of an infinite path \ if and only if \  $M$ is semi-Gorenstein-projective.

\vskip5pt

{\rm(1$''$)} \ $[M]$ is of $\mho$-type $-\Bbb N$ \ if and only if \ $M$ is semi-Gorenstein-projective, but
   not Gorenstein-projective.

\vskip5pt

{\rm(2)} \ $[M]$ is the end of a path of length $t\ge 1$ \ if and only if \
  $M$ is $t$-torsionfree  for $1\le i \le t$,  \ if and only if \  $\mho^{i-1} M$ is torsionless for $1\le i \le t$.

 In particular,  $[M]$ is the end of an arrow \ if and only if \  $M$ is torsionless$;$
 and  $[M]$ is the end of a path of length $2$ \ if and only if \  $M$ is reflexive.

\vskip5pt

{\rm(2$'$)} \ $[M]$ is the end of an infinite path \ if and only if \
  $M$ is $\infty$-torsionfree, \ if and only if \
 $M$ is reflexive and $M^*$ is semi-Gorenstein-projective.

\vskip5pt

{\rm(2$''$)} \ $[M]$ is of $\mho$-type $\Bbb N$ \ if and only if \ $M$ is $\infty$-torsionfree, but
   not Gorenstein-projective.

\vskip5pt

{\rm(3)} \  $[M]$ is the start of an infinite path and also the end of an infinite path
\ if and only if \  $M$ is Gorenstein-projective.

$M$ is of $\mho$-type $\Bbb Z$ \ if and only if \  $M$ is Gorenstein-projective
 and not $\Omega$-periodic.

$M$  is of $\mho$-type $\widetilde{\Bbb A}_n$ for some $n\ge 0$\ if and only if \
    $M$ is Gorenstein-projective and $\Omega$-periodic.
	
\vskip5pt

{\rm(4)} \ $A$-duality provides a bijection between the isomorphism classes of the
  reflexive indecomposable $A$-modules
  of $\mho$-type $\Bbb A_n$ and the isomorphism classes of the reflexive indecomposable $A^{\op}$-modules
  of $\mho$-type $\Bbb A_n$.

  Thus, for any $n\ge 3$,
$A$ has $\mho$-components of form $\Bbb A_n$ \ if and only if \
$A^{\op}$ has $\mho$-components of form $\Bbb A_n$.
	
\vskip5pt

{\rm(5)} \
 $A$-duality provides a bijection between the isomorphism classes of the
  reflexive indecomposable $A$-modules
  of $\mho$-type $\Bbb N$ and the isomorphism classes of the reflexive indecomposable $A^{\op}$-modules
  of $\mho$-type $-\Bbb N$.

Thus, $A$ has $\mho$-components of form $\Bbb N$ \ if and only if \
$A^{\op}$ has $\mho$-components of form $-\Bbb N$.\end{thm}

\begin{rem} \label{quiver2}  1. \ The $\mho$-quiver shows that an indecomposable module $M$ is
Gorenstein-projective if and only if both $M$ and $\Tr M$ are semi-Gorenstein-projective,
if and only if $M$ is reflexive and both $M$ and $M^*$ are semi-Gorenstein projective.
See (1$'$), (2$'$) and (3).

\vskip5pt

2. \  An indecomposable non-projective module $M$ is
semi-Gorenstein-projective if there is an infinite $\mho$-path
starting in $M$; and $M$ is $\infty$-torsionfree if
there is an infinite $\mho$-path ending in $M$.

\vskip5pt

 3. \  $A$ is left (respectively, right) weakly Gorenstein
if and only if there are no modules of $\mho$-type $-\Bbb N$ (respectively, $\Bbb N$), see (1$''$) (respectively,
(2$''$) and (5)).
\end{rem}	

\section{\bf Gorenstein-projective modules over the Nakayama algebras}

Claus Michael Ringel always emphasizes roles of examples in research.
When he began to study the Gorenstein-projective modules, he looked for more examples.
However, it seemed that not much was known about concrete construction and computation of this kind of modules at that time.
So he tried to first look at the Nakayama algebras $A$. With his high techniques of quivers,
C. M. Ringel could soon describe the structure of the category of
the Gorenstein-projective $A$-modules; and by using the so-called resolution quiver,
he provides a fast algorithm to obtain the Gorenstein-projective $A$-modules, and to judge whether $A$ is Gorenstein or CM-free.
Thus, this work becomes his first published paper on the Gorenstein-projective modules.

\vskip5pt

A module $M$ over an artin algebra is {\it uniserial},  if the set of submodules  is totally ordered by inclusion.
An artin algebra $A$ is Nakayama, if both the indecomposable projective modules and the indecomposable injective modules are uniserial.
The quiver of a Nakayama algebra is either $\mathbb{A}_n$ with linear order,
or a cycle with cyclic orientation. If it is $\mathbb{A}_n$,
then the Gorenstein-projective modules are precisely the projective modules. So, assume that  $A$ is a connected Nakayama algebra whose quiver is a cycle with cyclic orientation, or equivalently,
$A$ is a connected Nakayama algebra without a simple projective module.
Let $\gp_{0}(A)$ be the subcategory of all Gorenstein projective modules without indecomposable projective direct summands.

\vskip5pt

C. M. Ringel introduced  {\it the Gorenstein core} $\mathcal{C}$, which turns out to be essential in describing the structure of the category of
the Gorenstein-projective $A$-modules. By definition $\mathcal{C}$ is the full subcategory of $A$-mod  of
direct sums of indecomposable non-projective Gorenstein-projective modules and their projective covers. Then
$$\gp_{0}(A) \subseteq  \mathcal C \subseteq \gp(A).$$

\vskip5pt

Recall that for an artin algebra $A$, let $A\mbox{-}\underline{{\rm mod}}$ denote the stable category of $A\mbox{-}{\rm mod}$, i.e.,
the objects of $A\mbox{-}\underline{{\rm mod}}$ are precisely $A$-modules, and
$\Hom_{A\mbox{-}\underline{{\rm mod}}} (X,Y)$ is the quotient group of  $\Hom_{A\mbox{-}{\rm mod}}(X,Y)$ with respect to the subgroup
of those morphisms which factors through some projective $A$-modules.
Similarly, let $\underline{\gp}(A)$ denote the stable category of $\gp(A)$.

\vskip5pt

Usually $\gp(A)$ is not an abelian category.
Quite surprising, as the following result shows, after smartly replacing $\gp(A)$ by $\mathcal{C}$,
C. M. Ringel even could prove that $\mathcal{C}$ is equivalent to the module category of a self-injective Nakayama algebra.
In particular, $\mathcal{C}$ is an abelian Frobenius category, and the corresponding stable category is just $\underline{\gp}(A)$, up to a triangle-equivalence.

\begin{thm} \label{nakayama1} {\rm ([R1, Proposition 1])} \ Let $A$ be a connected Nakayama algebra without a simple projective module.
Then

 \vskip5pt

 $(1)$ \ The Gorenstein core $\mathcal{C}$ is an abelian subcategory of $A\mbox{-}{\rm mod}$ which is closed under extensions, projective covers and minimal left $A\mbox{-}$approximations.

\vskip5pt

 $(2)$ \  $\mathcal{C}$ is equivalent to $\Lambda\mbox{-}{\rm mod}$,  where $\Lambda$ is a self-injective Nakayama algebra and the inclusion functor $\mathcal{C}\hookrightarrow \gp(A)$ induces an equivalence
$\Lambda\mbox{-}\underline{{\rm mod}}\longrightarrow \underline{\gp}(A)$.

\vskip5pt

 $(3)$ \  If $\mathcal C\ne 0$ then $\Lambda$ is connected.
 \end{thm}

\vskip 5pt

Note that the inclusion functor $\mathcal{C}\hookrightarrow \gp(A)$ is exact, and the induced equivalence
$\Lambda\mbox{-}\underline{{\rm mod}}\longrightarrow \underline{\gp}(A)$ is a triangle-equivalence.

\vskip 5pt

\begin{thm} \label{nakayama2} {\rm ([R1, Proposition 2])} \ \ Let $A$ be a connected Nakayama algebra without a simple projective module,
and $\mathcal{E}\subseteq \gp_{0}(A)$ the set of non-zero modules such  that
no proper non-zero factor module is Gorenstein-projective. Then

\vskip5pt

$(1)$ \  $\mathcal{C}=\mathcal{F}(\mathcal{E})$, where $\mathcal{F}(\mathcal{E})$
the class of all modules with a filtration with factors in $\mathcal{E}$.
\vskip5pt

$(2)$ \ Assume that $\mathcal{C}\ne 0$.
Let  $E_{1}, \cdots, E_{g}$ be the representatives of the isomorphism classes in $\mathcal{E}$.
Then any simple module has multiplicity $1$ as a composition factor of $\bigoplus\limits_{i=1}^g E_i$. In particular, $E_{1}, \cdots, E_{g}$ are pairwise orthogonal bricks.

\vskip5pt

$(3)$ \ $\mathcal{E}$ is precisely the set of non-zero modules in $\gp_{0}(A)$ such that no proper non-zero submodule is Gorenstein projective.
\end{thm}

The elements in  $\mathcal E$ is precisely all the simple objects of $\mathcal C$. They are called {\it the elementary Gorenstein-projective modules}.
Note that $\mathcal C$ is the extension closure of $\gp_0(A).$

\vskip5pt

Note that the proof of  Theorem \ref{nakayama1} needs to use Theorem \ref{nakayama2}.

\vskip5pt

For a class $\mathcal M$ of modules, denote by $\add\mathcal M$ the smallest full subcategory of $A$-mod which contains $\mathcal M$ and is closed under direct sums and direct summands.

\begin{thm} \label{nakayama3} {\rm ([R1, Proposition 3])} \  Let $\mathcal{X}$ be the class of simple modules $S$  such that
the projective cover $P(S)\in \mathcal{C}$. Then $M\in \mathcal{C}$
if and only if ${\rm top}M\in {\rm add}\mathcal{X}$ and ${\rm soc}M\in {\rm add}(\tau^{-}\mathcal{X})$; and  if and only if ${\rm top}M\in {\rm add}\mathcal{X}$ and ${\rm top}\Omega M\in {\rm add}\mathcal{X}$.
\end{thm}

The Gorenstein core $\mathcal C$ may be obtained by deleting ray and corays from the Auslander - Reiten quiver
of $A$: the rays consisting of the indecomposable modules with top not
in $\mathcal X$, as well as the corays consisting of the indecomposable modules with socle not in $\tau^-\mathcal X$.

\vskip 5pt

An indecomposable projective module $P$ is {\it minimal}, if any proper non-zero submodule of $P$ is no longer projective, or equivalently, the projective dimension of
${\rm top} P$ is at least $2$.

\begin{thm} \label{nakayama4} {\rm ([R1, Proposition 4])} \  Let $A$ be a Nakayama algebra without a simple projective module, and $M$ an indecomposable non-projective $A\mbox{-}$module.
Then the following are equivalent.

\vskip5pt

$(1)$ \  $M$ is Gorenstein-projective.

\vskip5pt

$(2)$ \ All the projective modules occurring in a minimal projective resolution of $M$ are minimal.

\vskip5pt

$(3)$ \  There is an exact sequence $0\longrightarrow M\longrightarrow P_{n-1}\longrightarrow\cdots\longrightarrow P_{0}\longrightarrow M\longrightarrow 0$ such that all the modules $P_{i}$ are minimal projective.
\end{thm}

In fact, the exact sequence in Theorem \ref{nakayama4}(3) is a concatenation of $\mho$-sequences.

\vskip 5pt

For a Nakayama algebra $A$ without a simple projective module, the vertices of {\it the resolution quiver} $R = R(A)$ are the simple modules; and for every vertex $S$, there is an arrow from $S$ to $\tau {\rm soc} P(S)$.
The resolution quiver provides a visualization of the function from $S$ to $\tau {\rm soc} P(S)$, which has been studied by W. H. Gustafson [G].
C. M. Ringel realized the importance of this function, and named the quiver as  the resolution quiver.
Since any vertex in $R$ is the start of a unique arrow, in any component of $R$ there is a unique cycle with cyclic orientation.

\vskip 5pt

A vertex $S$ of $R$ is said to be {\it black} provided the projective dimension of $S$ is at least 2, or equivalently, $P(S)$ is a minimal
projective module; otherwise it is said to be {\it red}.

\vskip 5pt

Let $s$ be the number of simple $A$-modules, and $p$ the minimal length of the indecomposable projective $A$-modules.
Then $p\ge 2$, and $p > s$ if and only if no projective module is a brick.
As pointed out in [R1], in the following result one needs the condition $p>s$: for $p\le s$ there are some irregularities.

\vskip 5pt

An artin algebra $A$ is  {\it CM-free} ([Chen3]) if  all the Gorenstein projective modules are projective.

\begin{thm} \label{nakayama5} {\rm ([R1, Proposition 5])} \  Let $A$ be a Nakayama algebra without a simple projective module. Assume that $p> s$.

\vskip 5pt

$(1)$ \ $A$ is a Gorenstein algebra if and only if any cycle in $R$ consists of black vertices.

\vskip 5pt

$(2)$ \  \ $A$ is CM-free if and only if any cycle in $R$ contains at least one red vertex.
\end{thm}

\begin{exm} \label{nakayama6} {\rm ([R1, p. 243])} \  Let $Q$ be the quiver of type $\widetilde A_4$ with cyclic orientation. Instead of writing down the relations,
it suffices to mention the length of the indecomposable projectives. Here $|P(i)| = 13$ for $i=1,2$ and $|P(i)|= 12$ for $i=3,4,5$.  There are two elementary
Gorenstein-projective modules, namely $E(1)$ with composition factors $1,2,3$ and
$E(4)$ with composition factors $4,5$. The right picture below shows the
support of these modules $E(1)$ and $E(4)$, thus the $\Ext$-quiver for
$\mathcal  E = \{E(1),E(4)\}$ will be the quiver of type $\widetilde{\Bbb A}_1$ with
cyclic orientation.
$$
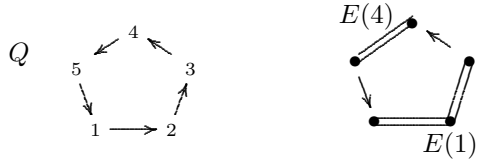
\hbox{\beginpicture
\setcoordinatesystem units <1cm,1cm>
%============================================
\put{\beginpicture
\put{$Q$} at -1 1
\put{$\ssize 1$} at 0 0
\put{$\ssize 5$} at -.25 .8
\put{$\ssize 4$} at 0.5 1.3
\put{$\ssize 3$} at 1.25 .8
\put{$\ssize 2$} at 1 0
\arr{0.2 0}{0.8 0}
\arr{-.2 .6}{-.05 .2}
\arr{1.05 .2}{1.2 .6}
\arr{0.3 1.2}{0 1}
\arr{1 1}{0.7 1.2}

\endpicture} at 0 0
%============================================
\put{\beginpicture
%\put{$Q'$} at -1.2 1
\put{$\bullet$} at 0 0
\put{$\bullet$} at -.25 .8
\put{$\bullet$} at 0.5 1.3
\put{$\bullet$} at 1.25 .8
\put{$\bullet$} at 1 0
\plot 0 0.05  1 0.05 /
\plot 0 -.04  1 -.04 /
\arr{-.2 .6}{-.05 .2}
\plot 1.05 0  1.3 .8 /
\plot 0.95 0  1.2 .8 /
\arr{1 1}{0.7 1.2}
\plot -.2 .8  0.5 1.3 /
\plot -.2 .9  0.5 1.4 /
\put{$E(1)$} at 1 -.3
\put{$E(4)$} at -.1 1.4

\endpicture} at 4 0
\endpicture}
$$

Draw twice the Auslander-Reiten
quiver of $A$. The left dashed boundary has to be identified with the right dashed boundary, and label the simple module $S(i)$ by $i$.
On the left, bullets are the indecomposable objects in $\mathcal  C$, and the encircled ones
are the elementary Gorenstein-projective modules $E(1)$ and $E(4)$.
In the middle, one shades the rays and the
corays which have to be deleted: note that
$\mathcal  X = \{S(1),S(4)\},$ thus $\tau^-\mathcal  X = \{S(5),S(3)\},$ this means that
one has to delete the corays ending in $S(2), S(3),S(5)$, and the rays starting in
$S(1), S(2), S(4)$. On the right, it is the Auslander-Reiten quiver of $\mathcal  C$.

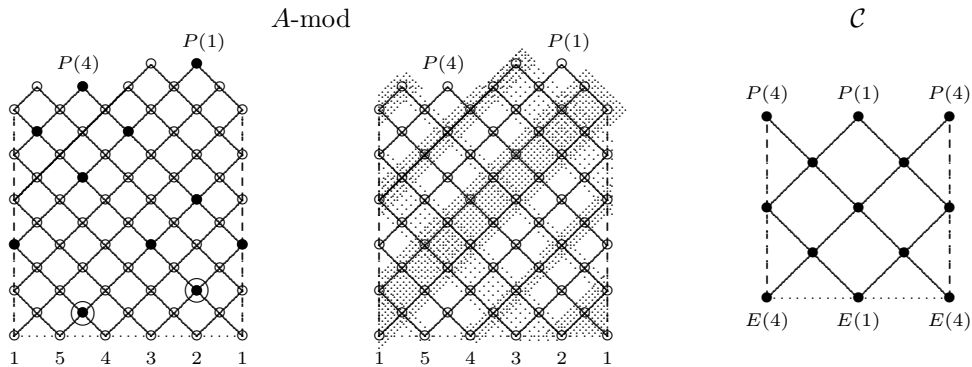
$$
\hbox{\beginpicture
\setcoordinatesystem units <.3cm,.3cm>
%============================================
\put{\beginpicture
\multiput{$\circ$} at
   0 0  0 2  0 4  0 6  0 8  0 10
   1 1  1 3  1 5  1 7  1 9  1 11
   2 0  2 2  2 4  2 6  2 8  2 10
   3 1  3 3  3 5  3 7  3 9  3 11
   4 0  4 2  4 4  4 6  4 8  4 10
   5 1  5 3  5 5  5 7  5 9  5 11
   6 0  6 2  6 4  6 6  6 8  6 10  6 12
   7 1  7 3  7 5  7 7  7 9  7 11
   8 0  8 2  8 4  8 6  8 8  8 10  8 12
   9 1  9 3  9 5  9 7  9 9  9 11
   10 0  10 2  10 4  10 6  10 8  10 10 /
\plot 0 0  10 10  8 12  0 4  4 0  10 6  5 11
   0 6  6 12  10 8  2 0  0 2  9 11 /
\plot  8 0   10 2  1 11  0 10  10 0  /
\plot 0 6  6 0  10 4  3 11  0 8  8 0 /
\setdashes <1mm>
\plot 0 0  0 10 /
\plot 10 0  10 10 /
\setdots <1mm>
\plot 0 0  10 0 /
\multiput{$\bullet$} at 0 4  1 9  3 1  3 7  3 11  5 9  6 4  8 2  8 6  8 12  10 4 /
\multiput{$\bigcirc$} at 3 1  8 2 /
\put{$\ssize {P(1)}$} at 8.3 13
\put{$\ssize {P(4)}$} at 2.8 12
\put{$\ssize 1$} at 10 -1
\put{$\ssize 2$} at 8 -1
\put{$\ssize 3$} at 6 -1
\put{$\ssize 4$} at 4 -1
\put{$\ssize 5$} at 2 -1
\put{$\ssize 1$} at 0 -1
\endpicture} at 0 0
%============================================
\put{\beginpicture
\multiput{$\circ$} at
   0 0  0 2  0 4  0 6  0 8  0 10
   1 1  1 3  1 5  1 7  1 9  1 11
   2 0  2 2  2 4  2 6  2 8  2 10
   3 1  3 3  3 5  3 7  3 9  3 11
   4 0  4 2  4 4  4 6  4 8  4 10
   5 1  5 3  5 5  5 7  5 9  5 11
   6 0  6 2  6 4  6 6  6 8  6 10  6 12
   7 1  7 3  7 5  7 7  7 9  7 11
   8 0  8 2  8 4  8 6  8 8  8 10  8 12
   9 1  9 3  9 5  9 7  9 9  9 11
   10 0  10 2  10 4  10 6  10 8  10 10 /
\plot 0 0  10 10  8 12  0 4  4 0  10 6  5 11
   0 6  6 12  10 8  2 0  0 2  9 11 /
\plot  8 0   10 2  1 11  0 10  10 0  /
\plot 0 6  6 0  10 4  3 11  0 8  8 0 /
\setdashes <1mm>
\plot 0 0  0 10 /
\plot 10 0  10 10 /
\setdots <1mm>
\plot 0 0  10 0 /
%\multiput{$\bullet$} at 0 4  1 9  3 1  3 7  3 11  5 9  6 4  8 2  8 6  8 12  10 4 /
\put{$\ssize {P(1)}$} at 8.3 13
\put{$\ssize {P(4)}$} at 2.8 12
\put{$\ssize 1$} at 10 -1
\put{$\ssize 2$} at 8 -1
\put{$\ssize 3$} at 6 -1
\put{$\ssize 4$} at 4 -1
\put{$\ssize 5$} at 2 -1
\put{$\ssize 1$} at 0 -1
\setshadegrid span <.4mm>
\vshade -.2  -.9  2.5  <,z,,> 9.15  8.5  11.9 <z,,,> 11 10  10 /
\vshade 7 0 0  10.2 0 3 /
\vshade -.2 9 10.7 <,z,,> 1.2 10.5  11.8 <z,,,> 2 11 11 /
\vshade -.2 5 6.7 <,z,,> 6.2 11.5  12.8 <z,,,> 7 12 12 /
\vshade 3.5 -.2 .2 <,z,,> 4.3 -.3 0.8  <z,,,> 10.3 5.6 6.9 /

\setshadegrid span <.6mm>
\vshade -.2 5.4 8.8 <,z,,> 5.9 -.3 2.7 <z,,,> 8.5 -.2 0 /
\vshade -.2 1.3 2.7 <,z,,> 1.9 -.3 0.7 <z,,,> 2.6 -.2 0.2 /
\vshade 1 10.2 11.7  10.2 1.3 2.3 /
\vshade 4.3 11 11 <,z,,> 6.2 9.2  12.7  <z,,,> 10.3 5.2  8.7 /
\endpicture} at 16 0
%============================================
\put{\beginpicture
\setcoordinatesystem units <.6cm,.6cm>
\multiput{$\bullet$} at
   0 0  2 0  4 0
   1 1  3 1
   0 2  2 2  4 2
   1 3  3 3
   0 4  2 4  4 4
   /
\plot 0 0  4 4  /
\plot 0 4  4 0 /
\plot 2 0  4 2  2 4  0 2  2 0 /
\setdashes <1mm>
\plot 0 0  0 4 /
\plot 4 0  4 4 /
\setdots <1mm>
\plot 0 0  4 0 /
\put{$\ssize {P(1)}$} at  2 4.5
\multiput{$\ssize {P(4)}$} at 0 4.5  4 4.5 /
\put{$\ssize {E(1)}$} at 2 -.5
\multiput{$\ssize {E(4)}$} at 0 -.5  4 -.5 /
\endpicture} at 32 -.4
\put{$A$-mod} at 8 8
\put{$\mathcal  C$} at 32 8

\endpicture}
$$

\vskip5pt

The resolution quiver $R$ has two components:
$$
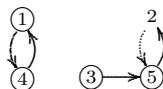
\hbox{\beginpicture
\setcoordinatesystem units <1cm,1cm>
%================================14solid
\put{\beginpicture
\setcoordinatesystem units <.8cm,.8cm>
\put{$\ssize 4$} at 1 0
\put{$\ssize 1$} at 1 1
\multiput{$\bigcirc$} at 1 0  1 1 /
\arr{0.89 0.22}{0.9 0.2}
\arr{1.11 0.78}{1.1 0.8}
\setquadratic
\plot 1.1 0.2 1.2 0.5 1.1 0.8 /
\plot 0.89 0.2 0.8 0.5 0.9 0.8 /

\endpicture} at 0 0

%================================25halfdotted
\put{\beginpicture
\setcoordinatesystem units <.8cm,.8cm>
\put{$\ssize 5$} at 1 0
\put{$\ssize 2$} at 1 1
\put{$\ssize 3$} at 0 0
\arr{0.2 0}{0.8 0}
\multiput{$\bigcirc$} at 1 0  0 0  /
\arr{0.89 0.22}{0.9 0.2}
\arr{1.11 0.78}{1.1 0.8}
\setquadratic
\plot 1.1 0.2 1.2 0.5 1.1 0.8 /
\setdots <.4mm>
\plot 0.89 0.2 0.8 0.5 0.9 0.8 /
\endpicture} at 1.3 0
\endpicture}
$$
The black vertices $1,3,4,5$ have been encircled. The arrow $2 \longrightarrow 5$ has been dotted in
order to stress that it starts at a red vertex. This Nakayama algebra is neither Gorenstein nor CM free. \end{exm}

\section{\bf Representations of quivers over the algebra of dual numbers}

Instead of the representations of quiver $Q$ over field $k$,
one can study the representations of $Q$ over algebra $A$, or equivalently, the $\Lambda$-modules, where $\Lambda = A\otimes_kkQ$.
Although the Gorenstein-projective $\Lambda$-modules can be described via the monomorphism category of the Gorenstein-projective $A$-modules ([LuoZ]),
it is still hard to determine all the indecomposable Gorenstein-projective $\Lambda$-modules.
C. M. Ringel suggests to look at the representations of $Q$ over $A$, where $A$ is the algebra of dual numbers.
Very surprising, in this case, using the covering theory developed by Gabriel school, C. M. Ringel could prove that
the indecomposable non-projective Gorenstein-projective $\Lambda$-modules are in one to one correspondence with the indecomposable $kQ$-modules.

\vskip5pt

Let $Q$ be a finite connected acyclic quiver, $A = k[\epsilon] = k[x]/\langle x^2\rangle$, and $\Lambda = kQ[\epsilon] = (kQ)[x]/\langle x^2\rangle$. Then $\Lambda = A\otimes_k kQ.$

\subsection{\bf Perfect differential $kQ$-modules}

Let $R$ be a ring. {\it A differential {\rm (}left{\rm)} $R$-module} is a pair $(N,\epsilon)$,
where $N$ is a (left) $R$-module and $\epsilon\in \End_{R}(N)$ with $\epsilon^2 = 0$.
For differential $R$-modules $(N,\epsilon)$ and $(N',\epsilon')$, a morphism $f: (N,\epsilon) \longrightarrow (N',\epsilon')$
is an $R$-module homomorphism  $f:N \longrightarrow N'$ such that $\epsilon' f = f\epsilon$.
Denoted by ${\operatorname{diff}}(R)$ the category of differential $R$-modules. Then
${\operatorname{diff}}(R)$  is just the category of left $(R[x]/\langle x^2\rangle)$-modules.

\vskip5pt

A differential $R$-module $(N,\epsilon)$ is {\it perfect} if $N$ is a finitely generated projective $R$-module.
Denote by ${\operatorname{diff}}_{\operatorname{perf}}(R)$ the full subcategory of ${\operatorname{diff}}(R)$ consisting of perfect differential $R$-modules.

\vskip5pt

A morphism $f: (N,\epsilon) \longrightarrow (N',\epsilon')$ of differential $R$-modules is {\it homotopic to zero} if there exists an $R$-module homomorphism
$h:N \longrightarrow N'$ such that $f = h\epsilon + \epsilon'h$.
Denote by  $\underline{\operatorname{diff}}_{\operatorname{perf}}(R)$ the corresponding homotopy category.
\vskip5pt

Consider homology functor $H: {\operatorname{diff}}(R)\longrightarrow R\mbox{-}{\rm Mod},   \ H(N,\epsilon)=\Ker \epsilon/\Im \epsilon.$
It vanishes on the morphisms which are homotopic to zero.
 If $R$ is left noetherian, then  $H: {\operatorname{diff}}_{\operatorname{perf}}(R)\longrightarrow R\mbox{-}{\rm mod}.$

\vskip5pt

Let $B$ and $C$ be finite-dimensional $k$-algebras. Then ([AR3, 2.2])
$${\rm max}\{{\rm inj.dim.}_BB, \ {\rm inj.dim.}_CC\} \le {\rm inj. dim.}(B\otimes_kC) \le {\rm inj.dim.}_BB + {\rm inj.dim.}_CC.$$
Thus $B\otimes_kC$ is Gorenstein (self-injective, respectively) if and only if $B$ and $C$ are Gorenstein (self-injective, respectively).
So, in our case
$\Lambda = k[x]/\langle x^2\rangle\otimes_k kQ$ is a Gorenstein algebra with ${\rm inj.dim.}_\Lambda\Lambda\le 1$. For short we call such an algebra $1$-Gorenstein.
For a $1$-Gorenstein algebra, Gorenstein-projective modules are exactly the torsionless modules.
An important observation is

\begin{lem}\label{perfect} \ {\rm(\cite[2.1]{RZ3})} \ Let $R$ be left noetherian and left hereditary.
A differential $R$-module $(N,\epsilon)$ is perfect if and only if $N$ is a finitely generated $R$-module, and $(N,\epsilon)$ is a torsionless $(R[x]/\langle x^2\rangle)$-module.\end{lem}

Let $\mathcal T$ be an additive category and $F: \mathcal T\longrightarrow \mathcal T$ an automorphism. By definition, {\it the orbit category} $\mathcal T/F$
has the same objects as $\mathcal T$ and
$$\Hom_{\mathcal T/F}(X, Y)=\bigoplus\limits_{n\in Z}\Hom_{\mathcal T}(X, F^n Y).$$
Note that even if $\mathcal T$ is a triangulated category and $F$ is a triangle functor, in general
$\mathcal T/F$ has no triangulated structure such that the canonical functor $\mathcal T\longrightarrow \mathcal T/F$ is a triangle functor ([Kel2]).
Denote by $D^b(kQ\mbox{-}\mod)$ the bounded derived category of (finite-dimensional) $kQ$-modules.
Then the orbit category $D^b(kQ\mbox{-}\mod)/[1]$  is triangulated, where $[1]$ is the shift functor.
This follows from a general criterion given by B. Keller [Kel2] (it can be also seen from the following result).

\vskip5pt

A module $M$ is  {\it strongly Gorenstein-projective} ([BM])
provided that there is a complete projective resolution
\ $\cdots \longrightarrow  P \stackrel f \longrightarrow P \stackrel f \longrightarrow  P\longrightarrow \cdots $ such that $M = \Ker f$.

\begin{thm}\label{RZ3.1} \ {\rm(\cite[Theorem 1, Proposition 4.12]{RZ3})} \ Let $k$ be a field, $Q$ a finite connected acyclic quiver, and
$\Lambda = kQ[x]/\langle x^2\rangle = k[x]/\langle x^2\rangle\otimes_k kQ.$
Then

\vskip5pt

$(1)$ \  The algebra $\Lambda$ is $1$-Gorenstein$;$ the category $\gp(\Lambda)$ of Gorenstein-projective $\Lambda$-modules is exactly the category ${\operatorname{diff}}_{\operatorname{perf}}(kQ)$
of perfect differential $kQ$-module. Moreover, every  Gorenstein-projective $\Lambda$-module is strongly Gorenstein-projective.

\vskip5pt

$(2)$ \   The stable category of Frobenius category ${\operatorname{diff}}_{\operatorname{perf}}(kQ)$ is precisely the homotopy category $\underline{\operatorname{diff}}_{\operatorname{perf}}(kQ)$,
which is also triangle equivalent to the orbit category $D^b(kQ\mbox{-}\mod)/[1]$.

\vskip5pt

$(3)$ \  The homology functor $H: {\operatorname{diff}}_{\operatorname{perf}}(kQ) \longrightarrow kQ\mbox{-}\mod$ is a full and dense functor which furnishes a bijection between the
indecomposables in the homotopy category $\underline{\operatorname{diff}}_{\operatorname{perf}}(kQ)$ and those in $kQ\mbox{-}\mod$.
It yields a quiver embedding $\iota$ of the
Auslander-Reiten quiver of $kQ\mbox{-}\mod$ into the Auslander-Reiten quiver of
the homotopy category $\underline{\operatorname{diff}}_{\operatorname{perf}}(kQ)$.
\end{thm}

\subsection{\bf Indecomposable perfect differential $kQ$-modules}
One can regard a $kQ\mbox{-}$module as a $(kQ[x]/\langle x^2\rangle)\mbox{-}$module annihilated by $x.$ Since $\Lambda = kQ[x]/\langle x^2\rangle$ is Gorenstein, it follows that
any $\gp(\Lambda)$ is contravariantly finite in $\Lambda$-mod ([EJ2, Theorem 11.5.1]). Thus, each $kQ\mbox{-}$module admits a minimal right $\gp(\Lambda)$-approximation.

\vskip5pt

An algebra is {\it CM-finite} if it has only finitely many indecomposable Gorenstein-projective modules.
The following result shows a remarkable one to one correspondence between the indecomposable non-projective Gorenstein-projective $\Lambda$-modules  and
the indecomposable representations of $Q$, via the homology functor $H$ and the inverse construction $\eta$.

\begin{thm}\label{RZ3.2} \ {\rm(\cite[Theorem 2]{RZ3})} \ Let $k$, $Q$, $\Lambda$ be as in {\rm Theorem \ref{RZ3.1}}, and
$H: {\operatorname{diff}}_{\operatorname{perf}}(kQ) \longrightarrow kQ\mbox{-}\mod$ the homology functor. Then

\vskip5pt

$(1)$ \ For $(N, \epsilon) \in {\operatorname{diff}}_{\operatorname{perf}}(kQ)$, then $H(N, \epsilon) = 0$ if and only if $(N, \epsilon)$ is a projective $\Lambda$-module$;$ if and only if $(N, \epsilon)\cong (P\oplus P, \left(\begin{smallmatrix}0 & 1 \\ 0 & 0\end{smallmatrix}\right))$
for some projective $kQ$-module $P$.

\vskip5pt

$(2)$ \ The homology functor $H: {\operatorname{diff}}_{\operatorname{perf}}(kQ) \longrightarrow kQ\mbox{-}\mod$
induces a bijection between the indecomposable non-projective Gorenstein-projective $\Lambda$-modules
and the indecomposable $kQ$-modules.
The inverse bijection $\eta$ is given by taking the minimal right $\gp(\Lambda)$-approximation of an indecomposable $kQ$-module.
In particular, $\Lambda$ is CM-finite if and only if $Q$ is of Dynkin type.

\end{thm}

The map $\eta$  is not functorial, but stably functorial.
The $\Lambda$-maps $f$ such that $H(f) = 0$ are quite complicated. They were called {\it the ghost maps}, and discussed in [RZ3, Section 6].

\section{\bf From submodule categories to preprojective algebras}

\subsection{\bf Preprojective algebras of type $\mathbb A_n$}

C. M. Ringel has smartly observed that the module category of the preprojective algebra $\Pi_{n-1}$ is a factor category of the
submodule category $\mathcal{S}(n)$, while $\mathcal{S}(n)$ is the category of the Gorenstein-projective $T_2(\Lambda_n)$-modules ([LiZ]).
This process produces an abelian category from a factor  category of a triangulated category.

\vskip5pt

Let $k$ be a field, $n$ an integer with $n\ge 2$. and $\mathcal{S}(n)$ the category of invariant subspaces of nilpotent operators with nilpotency index at most $n$,
i.e., an object of $\mathcal{S}(n)$ is a triple $(V, U, T)$, where $U$ is a subspace of finite-dimensional $k$-space $V$ and $T\in \End_k(V)$,
such that $T^n = 0$ and $T(U)\subseteq U$; and a morphism from $(V,U,T)$ to $(V',U',T')$ is a $k$-linear map $f: V\longrightarrow V'$ such that $f(U)\subseteq U'$ and $fT=T'f$.
See C. M. Ringel and M. Schmidmeier [RS1].

\vskip5pt

Let $\Lambda_n = k[x]/\langle x^n\rangle$. For $0\leq i\leq n$, denote by $[i]$ the  $\Lambda_n$-module $\langle x^{n-i}\rangle/\langle x^n\rangle$. Then
$[i], \ 1\le i\le n,$ give all the indecomposable  $\Lambda_n$-modules. Put $E = \bigoplus\limits_{i=1}^n[i]$.

\vskip5pt

Let $A_n$ be the Auslander algebra of $\Lambda_n$, i.e.,  $A_n=\End_{\Lambda_n}(E)^{\rm op}$.
Then ${\rm gl.dim} A_n = 2 = {\rm dom.dim} A_n$ ([ARS, VI.5]); and $A_n$ is representation-finite if and only if $n\le 3$,
and there are only finitely many indecomposable torsionless $A_n$-modules if and only if $n\le 5$ ([DR, 7.1, 7.2]).
The algebra $A_n$ has a unique indecomposable projective-injective module $P(n) = \Hom_{\Lambda_n}(E, [n]) = \Hom_{\Lambda_n}(E, \Lambda_n)$ ([ARS, VI, 5]), and $A_n$ is quasi-hereditary.

\vskip5pt

Let $T_2(\Lambda_n)=\left(\begin{smallmatrix}\Lambda_n &\Lambda_n\\0 &\Lambda_n \end{smallmatrix}\right)$ be the upper triangular matrix algebra.
It is a $1$-Gorenstein algebra of infinite global dimension, and $T_2(\Lambda_n)\mbox{-}\mod$ is equivalent to the morphism category of $\Lambda_n\mbox{-}\mod$.
So, each $T_2(\Lambda_n)$-module will be identified with a triple $(X, \ Y, \ f)$, where $f: X\longrightarrow Y$ is a $\Lambda_n\mbox{-}$homomorphism;
and $\mathcal{S}(n)$ is precisely the full subcategory of $T_2(\Lambda_n)\mbox{-}\mod$ of triples $(X, \ Y, \ f)$, where $f: X\longrightarrow Y$ is a $\Lambda_n\mbox{-}$monomorphism.
Then $$\mathcal{S}(n) = \gp(T_2(\Lambda_n))$$ i.e.,  $\mathcal{S}(n)$ is precisely the category of Gorenstein-projective $T_2(\Lambda_n)$-modules ([LiZ, 1.1]).
Note that $T_2(\Lambda_n)$ is representation-finite  if and only if $A_n$ is representation-finite  ([ARS, VI, 5]). So,
$T_2(\Lambda_n)$ is representation-finite if and only if $n\le 3$. Moreover, if $n =4$ or $5$, then  $T_2(\Lambda_n)$ is representation-infinite and CM-finite ([LiZ, 1.3]).

\vskip5pt

Let $Q = (Q_0, Q_1)$ be a finite quiver, and $\overline{Q}$ the quiver obtained by adding an inverse arrow $\alpha^*$ for each $\alpha\in Q_1$.
Denote by $k \overline{Q}$ the path algebra of $\overline{Q},$ and $\langle\sum\limits_{\alpha \in Q_1}(\alpha \alpha^*-\alpha^*\alpha) \rangle$ the ideal of $k \overline{Q}$ generated by
$\sum\limits_{\alpha \in Q_1}[\alpha^*, \alpha]$, where $[\alpha^*, \alpha] = \alpha^*\alpha - \alpha\alpha^*$.
By definition the preprojective algebra of the quiver $Q$ is
$$\Pi(Q) =k \overline{Q}/\langle\sum\limits_{\alpha \in Q_1}[\alpha^*, \alpha]\rangle.$$ (If there are no arrows in $Q$ then $\Pi(Q) = k^{|Q_0|}$.)
Then $\Pi(Q)$ is finite-dimensional if and only if $Q$ is a disjoint union of Dynkin diagrams, and in this case $\Pi(Q)$ is self-injective ([DR, 6.1, 6.2]).
An alternative characterization of $\Pi(Q)$ is given by C. M. Ringel [R4, Theorem A]:
it is isomorphic to the tensor algebra of $kQ$-$kQ$-bimodule $\Ext_{kQ}^1 (D(kQ_{kQ}), \  _{kQ}kQ)$, where $D(kQ) = \Hom_k(kQ, k)$.

\vskip5pt

If $Q$ is the quiver of type $\mathbb{A}_n$ with linear order, then  $\Pi(Q)$ will be denoted by
$\Pi_n$ (thus $\Pi_1 = k$). An important relation between the Auslander algebra $A_n$ and  the preprojective algebras $\Pi_{n-1}$ for $n\ge 2$ is as follows.
Let $e_n$ be the idempotent of $A_n$ such that $P(n) = A_ne_n$. Then
$$\Pi_{n-1}\cong A_n/\langle e_n\rangle$$
See V. Dlab and C. M. Ringel [DR, Theorems 3, 4, and Section 7]. This induces an functor
$$\delta: A_n\mbox{-}\mod \longrightarrow \Pi_{n-1}\mbox{-}\mod,  \ \ M\mapsto M/A_neM.$$
Note that $\Pi_{n}$ is representation-finite  if and only if $n\le 4$ ([DR, 6.3]).

\vskip5pt

Define  two functors $F, G: \mathcal{S}(n)\longrightarrow \Pi_{n-1}\mbox{-}\mod$ as follows.
There are two functors $\iota, \epsilon: \mathcal{S}(n)\longrightarrow T_2(\Lambda_n)\mbox{-}\mod$,
where $\iota$ is the natural embedding, and $\epsilon$ sends each $(X,Y, f)\in \mathcal{S}(n)$ to $(Y, \ Y/f(X), \ \pi)\in T_2(\Lambda_n)\mbox{-}\mod$,
where $\pi$ is the canonical projection. Also, there is a functor
$$\alpha=\Coker \Hom_{\Lambda_n}(E,-): \ T_2(\Lambda_n)\mbox{-}\mod\longrightarrow A_n\mbox{-}\mod, \ \
(X,Y, f)\mapsto \Coker \Hom_{\Lambda_n}(E, f).$$
By definition $F=\delta\alpha\iota$ and $G=\delta\alpha\epsilon$. See the following diagram.
$$
%====================================
\hbox{\beginpicture
\setcoordinatesystem units <1.5cm,.7cm>
\put{$\mathcal{S}(n)$} at 0 0
\put{$T_2(\Lambda_n\mbox{-}\mod)$} at 2 0
\put{$A_n\mbox{-}\mod$} at 4 0
\put{$\Pi_{n-1}\mbox{-}\mod$} at 6 0
\arr{0.4 0.2}{1.3 0.2}
\arr{0.4 -.2}{1.3 -.2}
\arr{2.7 0}{3.5 0}
\arr{4.5 0}{5.4 0}
\put{$\iota$\strut} at .8 0.5
\put{$\epsilon$\strut} at .8 -.5
\put{$\alpha$\strut} at 3.1 0.3
\put{$\delta$\strut} at 4.9 0.3
\endpicture}
$$

\vskip5pt

For full subcategory $\mathcal B$ of an additive category $\mathcal A$, let $\mathcal A/\mathcal B$ denote the factor category of $\mathcal A$ with respect to $\mathcal B$, i.e.,
the objects of $\mathcal A/\mathcal B$ are the objects of $\mathcal A$, and
$\Hom_{\mathcal A/\mathcal B} (X,Y)$ is the quotient group of $\Hom_{\mathcal A}(X,Y)$ with respect to the subgroup
of those morphisms which factors through some object in $\mathcal B.$
If $\mathcal B$ is an additive category then so is $\mathcal A/\mathcal B$.

\vskip5pt

Let $F: \mathcal A\longrightarrow \mathcal B$ be an additive functor between additive categories.
Denote by $\Ker F$ the full subcategory of $\mathcal A$ of objects $X$ with $F(X)=0$. By definition, $F$ is  {\it objective}, if any morphism $f: \ X
 \longrightarrow Y$ with $F(f) = 0$ factors through an object $K\in \Ker F$.

\vskip5pt

\begin{fact}
Let $F: \mathcal A\longrightarrow \mathcal B$ be an additive functor between additive categories. If $F$ is full, dense and objective,
then $F$ induces an equivalence $\mathcal A/\Ker F\cong \mathcal B$ of categories.
\end{fact}

It turns out that  $\Pi_{n-1}$-mod is a factor category of the category of
the Gorenstein-projective $T_2(\Lambda_n)$-modules, in two ways.

\begin{thm}\label{RZ1.1} {\rm (\cite[Theorem 1]{RZ1})} \
$(1)$ \ Both the functors $F,G: \mathcal{S}(n)\longrightarrow \Pi_{n-1}\mbox{-}\mod$ are full, dense and objective.

\vskip5pt

$(2)$ \ The functor $F$ induces an equivalence $\mathcal{S}(n)/\Ker F\cong \Pi_{n-1}\mbox{-}\mod$, and $G$ induces an equivalence $\mathcal{S}(n)/\Ker G\cong \Pi_{n-1}\mbox{-}\mod$.
An indecomposable module in $\Ker F$ is exactly either $([i],[i])$ with $1\le i\le n$ or $([i], [n])$ with $0\le i\le n-1;$ an indecomposable module in $\Ker G$ is exactly either $([i],[i])$ with $1\le i\le n$ or $([0], [i])$ with $1\le i\le n$.

\vskip5pt

In particular, both $\mathcal{S}(n)/\Ker F$ and $\mathcal{S}(n)/\Ker G$ are abelian categories with enough projective objects. The indecomposable projective objects in $\mathcal{S}(n)/\Ker F$
are $([0],[i])$ with $1\le i\le n-1;$ and the indecomposable projective objects in $\mathcal{S}(n)/\Ker G$ are $([i], [n])$ with $1\le i\le n-1$.
\end{thm}

\begin{rem} \ Since $\mathcal{S}(n)\cong \gp(T_2(\Lambda_n)$, the stable category $\underline {\mathcal{S}(n)}$ is triangulated.
The indecomposable projective objects $([0], [n])$ and $([n],[n])$ of $\mathcal{S}(n)$ are in $\Ker F$ and $\Ker G$, so
$$\mathcal{S}(n)/\Ker F = \underline{\mathcal{S}(n)}/\Ker F, \ \ \ \mathcal{S}(n)/\Ker G =\underline{\mathcal{S}(n)}/\Ker G.$$
Thus the factor categories $\underline{\mathcal{S}(n)}/\Ker F$ and $\underline{\mathcal{S}(n)}/\Ker G$ of triangulated category $\underline {\mathcal{S}(n)}$
become abelian category $\Pi_{n-1}\mbox{-}\mod$. This gives examples that an abelian category is realized as a factor category of a triangulated category. Note that a famous theorem of
A. A. Beilinson - J. Bernstein - P. Deligne claims the heart of a $t$-structure of a triangulated category is an abelian category.
\end{rem}

\begin{thm} {\rm (\cite[Theorem 2]{RZ1})} \ Let $\pi: \Pi_{n-1}\mbox{-}\mod\longrightarrow \Pi_{n-1}\mbox{-}\underline{\mod}$ be the canonical projection functor, and $\Omega: \Pi_{n-1}\mbox{-}\underline{\mod}\longrightarrow \Pi_{n-1}\mbox{-}\underline{\mod}$ be the syzygy functor. Then
$$\pi F=\Omega \pi G: \ \mathcal S(n) \longrightarrow \Pi_{n-1}\mbox{-}\underline{\mod}.$$
\end{thm}

\subsection{\bf Objective triangle functors} \ All functors considered below are assumed to be covariant and additive between additive
categories. C. M. Ringel defined an objective functor. As we have seen, objective functors are used in the connection between submodule categories and preprojective algebras.
Recall that a functor  $F:  \mathcal A \longrightarrow \mathcal B$ is  {\it objective} if any morphism $f: \ X
 \longrightarrow Y$ with $F(f) = 0$ factors through
 an object $K$ with $F(K) = 0$; and $F$ is {\it sincere},
if $F$ sends non-zero objects to non-zero objects.
A functor is faithful if and only if it is objective and sincere.
The composition of objective functors is not necessarily objective;
while the composition of full, dense, objective functors is full, dense, and objective ([RZ1, Appendix]).

\vskip5pt

Triangle functors behave
quite different from general additive functors, and also from exact functors between abelian categories.
For examples, a full functor may not be
objective ([RZ1, Appendix]), but a full triangle
functor is objective ([RZ2, 4.4]). An exact functor is objective, whereas there
are sincere triangle functors which
are not objective ([RZ2, Section 8]). Also, for an exact functor $F$, $F$ is faithful if and only if $F$ is sincere;
for triangle functor $F$, $F$ is sincere and full imply that $F$ is faithful,  but $F$ is sincere does not imply that $F$ is faithful.

\vskip5pt

The following conditions for a triangle functor were proposed in [RZ2]:

\vskip5pt

${\rm (WSM)}$ \  For each morphism $u: X \longrightarrow Y$ such that $F(u)$ is a splitting monomorphism, there exists a morphism $u': Y \longrightarrow X'$
such that $F(u'u)$ is an isomorphism.

\vskip5pt

${\rm (I)}$ \ For each morphism $u: X \longrightarrow Y$ such that $F(u)$ is an isomorphism, there exists a morphism $u': Y \longrightarrow X$ such that
$F(u)^{-1} = F(u').$

\vskip5pt

${\rm (SM)}$ \  For each morphism $u: X \longrightarrow Y$ such that $F(u)$ is a splitting monomorphism, there exists a morphism $u': Y \longrightarrow X$ such
that $F(u'u) = {\rm Id}_{F(X)}$.

\vskip5pt

A functor $F$ satisfies both conditions (WSM) and (I) if and only if it satisfies condition ${\rm (SM)}$.

If a triangle functor $F$ is faithful or full,
then $F$ satisfies  conditions ${\rm (WSM)}$ and ${\rm (I)}$.

\vskip5pt

If $F: \mathcal
A\longrightarrow \mathcal B$ is a triangle functor, then ${\Ker}F$ is a
triangulated subcategory of $\mathcal A$, and one has the
Verdier quotient functor $V_F: \mathcal A \longrightarrow \mathcal A/\Ker F$. Here $\mathcal A/\Ker F$ is not the usual additive factor category,
instead, by definition, the triangulated category $\mathcal A/\Ker F$ is the localization of
$\mathcal A$ with respect to the compatible multiplicative system determined by $\Ker F$.
Since
$F({\Ker}F) = 0$, the universal property of the Verdier quotient
functor asserts that there exists a unique triangle functor
$\widetilde{F}: \mathcal A/{\Ker} F \longrightarrow \mathcal B$,
such that $F = \widetilde FV_F$. The functor $\widetilde F$ is
always sincere.

\begin{thm} \label{1.1} \  {\rm(\cite[1.1]{RZ2})} \ Let $F: \mathcal
A \longrightarrow \mathcal B$ be a triangle functor between
triangulated categories. Then the following are equivalent$:$

$\rm (1)$ \ \ $F$ is objective{\rm;}

$\rm (2)$ \ \ $F$ satisfies the condition ${\rm (WSM)};$

$\rm (3)$ \ \ The induced functor $\widetilde{F}: \mathcal
A/{\Ker}F \longrightarrow \mathcal B$ is faithful.
\end{thm}

In particular, if $\mathcal K$ is a triangulated subcategory of $\mathcal A$, then the Verdier
quotient functor $\mathcal A \longrightarrow \mathcal A/\mathcal K$ is objective (see also [Krau, 4.6.2]).

\vskip5pt

An additive category $\mathcal A$ is {\it Fitting} if for any morphism $a: X \longrightarrow X$
there is a direct decomposition $X = X'\oplus X''$ with $a(X')
\subseteq X', a(X'') \subseteq X''$ such that $a|_{X'}$
is an automorphism and $a|_{X''}$ is
nilpotent. For example, if $\mathcal A$ is a $\Hom$-finite
$k$-category, where $k$ is a field, and any object of $\mathcal A$
is a finite direct sum of objects with local endomorphism rings,
then $\mathcal A$ is Fitting.

\begin{thm} \label{1.2}  \  {\rm(\cite[1.2]{RZ2})} \ Let $F: \mathcal A \longrightarrow \mathcal B$ be a
triangle functor.

$(1)$ \ \ If the
Verdier quotient functor $V_F$ is full, then $F$ satisfies the condition ${\rm
(I)}$.

$(2)$ \ \  Assume that $F$ is objective or that $\mathcal A$ is a
Fitting category. Then  $V_F: \mathcal A \longrightarrow \mathcal A/\Ker F$ is full if and only if $F$ satisfies
${\rm (I)}$.
\end{thm}

Denote by $\ker(F)$ the class of morphisms $f$
with $F(f) = 0.$  Then $\ker(F)$ is an ideal
of $\mathcal A$, whereas $\Ker(F)$ is a subcategory. Consider the factor
category $\mathcal A/\ker(F)$ and the canonical projection functor $\pi_F: \mathcal A \longrightarrow
\mathcal A/\ker(F)$.

\begin{thm} \label{1.3}  \  {\rm(\cite[1.3]{RZ2})} \  Let $F: \mathcal A \longrightarrow \mathcal B$ be a
triangle functor between triangulated categories. Then the following
are equivalent$:$

$\rm (1)$ \ \  $F$ satisfies the condition ${\rm (SM)};$

$\rm (2)$ \ \ $F$ is objective and the
Verdier quotient functor $V_F: \mathcal A \longrightarrow \mathcal A/\Ker F$ is full{\rm;}

$\rm (3)$ \ \ There is an equivalence $\Phi:\mathcal A/\ker(F) \to \mathcal A/\Ker(F)$ such that $V_F
=\Phi \pi_F;$

$\rm (4)$ \ \ $F = F_2F_1$ where $F_1$ is a
full and dense triangle functor and $F_2$ is a faithful triangle
functor.
\end{thm}

\section{\bf Modules over short local algebras}

Short local algebras have a lot interest in commutative algebra.
C. M. Ringel stresses their special role in the representation theory of (non-commutative) algebras.

\vskip5pt

In this section $A$ is a finite-dimensional local algebra over field $k$ with radical $J$ and simple module $S$. For simplicity, we assume $A/J = k$.
A local algebra $A$ is {\it short} if $J^3 = 0.$
Denote by $|M|$ the length of the module $M$. Write $t(M) = |\top M|$. For a short local algebra $A$, put
$$e = |J/J^2|, \ \ a = |J^2|$$
and call $(e, a)$ the {\it Hilbert-type} of $A$. An $A$-module has Loewy length at most 2 if and only if it is annihilated by $J^2$.
In this case, put
$\bdim M = (t(M),|JM|)$.

\vskip5pt

\subsection{Existence of non-projective reflexive modules}

A non-zero module $M$ of Loewy length at most 2 is {\it solid} if
any endomorphism of $M$ is a scalar multiplication on $\soc M$ (thus, any non-invertible endomorphism vanishes on the socle). A solid module
is indecomposable. The following result gives a numerical condition for which non self-injective short local algebras $A$ should satisfy
if there is a non-projective reflexive $A$-module.

\begin{thm} \label{(e,a)reflexive} {\rm ([RZ8, Theorem 1.1, Corollary 4.4])} \  Let $A$ be a short local algebra which is not self-injective.
If there exists a non-projective reflexive module, then  $2\leq a\leq e-1$, $J_A$ and $_AJ$ are solid, and $J^{2}={\rm soc}{_{A}A}={\rm soc}A_{A}$.
\end{thm}

The bound $a \le e-1$ cannot be improved: If $1\le a \le e-1$, there exists an algebra of
Hilbert type $(e,a)$ with non-projective reflexive modules ([RZ8, 15.1]).
In general, ${}_AJ$ may be solid, whereas
$J_A$ is not solid ([RZ8, 4.8]).
	
\vskip5pt

Theorem \ref{(e,a)reflexive} assumes that there exists a non-projective reflexive module, thus
an $\mho$-path of length 2,
or equivalently, a non-projective module $M$ with $\Ext^i_A(M,A) = 0,$ for $1\le i \le 2.$
The next theorem shows that
the existence of a non-projective module $M$ with
$\Ext^i_A(M,A) = 0,$ for $1\le i \le 4$, yields a stronger assertion.

\begin{thm} \label{reflexive} {\rm ([RZ8, Theorem 1.2])} \
Let $A$ be a short local algebra which is not self-injective. Assume that  $M$ is an indecomposable,
reflexive and non-projective module with ${\rm Ext}^{i}_A(M, A) = 0$ for $i = 1, 2$. Then $2\leq a=e-1$,
and  ${\rm \bf dim}X= (t(M), at(M))$ for $X\in \{\Omega^{2}(M), \Omega(M), M, \mho(M)\}$.
\end{thm}

\subsection{\bf Short local algebras of dimension 6 with non-projective reflexive modules}
If there exists a non-projective reflexive module,  then $2\leq a\leq e-1$ (Theorem \ref{(e,a)reflexive}).
In this case, ${\rm dim}_k A \ge 6$; and if ${\rm dim}_k A = 6$ then the Hilbert type is $(3, 2)$.
In [R3] C. M. Ringel studies in detail such short local algebras, using his unpublished work [R5] on $K(3)$-modules.

\vskip5pt

A module is {\it uniform} if its socle is simple. Dually, a module is {\it local} if its top is simple.

\vskip 5pt

\begin{thm} \label{6(32)} {\rm ([R3, Theorem 1.1])} \   Let $A$ be a short local algebra over algebraically closed field $k$. The following conditions are equivalent:

\vskip 5pt

$(1)$ \ $J^{2}={\rm soc}{_{A}A}={\rm soc}A_{A}$ and there is no uniform left ideal with length $3$.

\vskip 5pt

$(2)$  \ There is a reflexive local module of length $3$ and Loewy length $2$.

\vskip 5pt

$(3)$  \ There is a non-projective reflexive module.
\end{thm}

C. M. Ringel also presents a 6-dimensional local algebra $A$ which is not self-injective and not short, such that there exists a non-projective reflexive $A$-module.

\vskip 5pt

\begin{exm} \label{simplereflnotselfinj} \ {\rm([R3, 12.1])} \ Let $A = k\langle x, y\rangle /\langle x^2, y^2, yxy\rangle.$ Then $M = Ax = A/Ax$ is a Gorenstein-projective module with $\Omega$-period 1. \end{exm}

\subsection{Existence of exact complex of projective modules}
A module $M$ is {\it bipartite} if $\soc M = J M$.
A module has Loewy length at most $2$ if and only if it is the direct sum of
a bipartite and a semisimple module.

\begin{thm} \label{exactcomplex} {\rm ([RZ8, Theorem 1.3])} \   Let $A$ be a short local algebra which is not self-injective, with a non-zero minimal
exact complex $P^{\bullet}= (P_{i}, d_{i})_{i}$ of projective modules. Then $1\leq a=e-1$.
Also, $M_{i}={\rm Im}d_{i}$ is bipartite for $i\ll 0$. Let $t_{i}=t(P_{i})=t(M_{i})$. There is $v\in \mathbb{Z}$ such that for $i\leq v$,
we have $t_{i}=t$ and ${\rm \bf dim} M_{i}=(t, at)$. There are just two possibilities:

\vskip5pt

{\rm Type I.} \ For all $i\in \Bbb Z,$ the module
 $M_i$ is bipartite with $\bdim M_i = (t,at)$ (thus $t_i = t$).

\vskip5pt

{\rm Type II.} \ One can choose $v$ in such a way that first, $t_{i+1} > t_i$ for $i\ge v$; second, the module $M_{v+1}$ is not bipartite; and third, $|JM_i| < at_i$ for $i > v$.
\end{thm}

\vskip5pt

For commutative rings Theorem \ref{exactcomplex} is
due to L. W. Christensen and O. Veliche [CV], where the case $a = 1$ does not occur.
But in general, the case $a = 1$ is possible ([RZ8, 9.2]).

\vskip5pt

 If $A$ is commutative, then the existence of a non-zero minimal exact
complex $P^\bullet$ of projective modules
implies that $J^2 = \soc A$; but in general it
does neither imply $J^2 = \soc {}_AA$, nor
$J^2 = \soc A_A$ ([RZ8, 9.2, 9.3]). If $J^2 = \soc {}_AA$,
then all $M_i$ with $i\le v$ are bipartite ([RZ8, 13.3]).

\subsection {\bf Semi-Gorenstein-projective and $\infty$-torsionfree modules.}

Theorems \ref{(e,a)reflexive} and \ref{exactcomplex} imply that if $A$ is a short local algebra which is not
self-injective, with a Gorenstein-projective module which is not
projective, then $2 \le a = e-1$.  The following result is a strengthening.

\begin{thm} \label{sgportrsgp} {\rm ([RZ8, Theorem 1.4])} \   Let $A$ be a short local algebra which is not self-injective. Assume that there exists
a non-projective indecomposable module  $M$ which is semi-Gorenstein-projective or $\infty$-torsionfree.
Then $2\leq a=e-1$ and $J^{2}= {\rm soc}_{A}A={\rm soc}A_{A}$. Moreover, let $t=t(M)$. One has in addition:

\vskip 5pt

$(1)$ \  If $M$ is torsionless and semi-Gorenstein-projective, then ${\rm \bf dim}\Omega^{i}(M)=(t, at)$ for all $i\gg 0$.

\vskip 5pt

$(2)$ \  If $M$ is $\infty$-torsionfree, then ${\rm \bf dim} \mho(M)=(t, at)$ for all $i\geq 0$.

\vskip 5pt

$(3)$ \ If $M$ is reflexive and semi-Gorenstein-projective, or if $M$  is $\infty$-torsionfree, then also ${\rm \bf dim}M^{*}=(t, at)$.

\vskip 5pt

$(4)$ \ If $M$ is Gorenstein-projective, then ${\rm \bf dim}X=(t, at)$ for $X\in \{\Omega^{i}(M), \ \mho^{i}(M), M^* \ | \ i\geq 0\}$.
\end{thm}

If $M$ is semi-Gorenstein-projective, but not reflexive, its Loewy length
may be 3; and if it is 2, one may have $\bdim M^* \neq \bdim M$  ([RZ8, 9.5]).

\subsection {\bf The Auslander-Reiten conjecture.}

For artin algebra $A$,  Auslander-Reiten ([AR1]) conjecture that if $M$ is a non-projective
semi-Gorenstein-projective module, then $\Ext^i_A(M,M)\neq 0$
for some $i\ge 1$.  For short local algebras this conjecture holds true in a stronger form.

\begin{thm} \label{ARconjecture} {\rm ([RZ8, Theorem 1.5])} \  Let $A$ be a short local algebra and $M$ a non-projective semi-Gorenstein-projective
module. Then ${\rm Ext}^{1}_A(M, M)\neq 0$. Moreover, if $A$ is not self-injective, then ${\rm Ext}^{i}_A(M, M)\neq 0$ for all $i\geq 1$.
\end{thm}

If $A$ is self-injective, then Theorem \ref{ARconjecture} is due to M. Hoshino [Hos].
For commutative short local rings, the Auslander-Reiten conjecture was proved by Huneke - \c Sega - Vraciu [HSV].
For a commutative self-injective short local algebra $A$, and a non-projective module $M$,
one even has $\Ext^i_A(M,M) \neq 0$ for  all $i\ge 1$ ([HSV]).
However, for $A$ non-commutative, this is not true: one
may have $\Ext^i_A(M, M) = 0$ for some or even for all $i \ge 2$ ([RZ8, A.19]).

\subsection{The $\Omega$-growth of modules} For a set $E\ne \emptyset$ of real numbers,
the least upper bound of $E$ and the greatest lower bound of $E$ exist, which are denoted by ${\rm sup}E$ and ${\rm inf}E$, respectively.

\vskip5pt

Let $\{a_n\}_{n\ge 1}$ be a sequence of real numbers. Put
$$b_l = {\rm sup}\{a_l, a_{l+1}, a_{l+2}, \cdots\}, \ \ l = 1, 2, 3, \cdots;  \ \ \ \mbox{and} \ \ \beta = {\rm inf}\{b_1, b_2, b_3, \cdots \}.$$
The number $\beta$ is called {\it the upper limits} of  $\{a_n\}_{n\ge 1}$, and write
\ $\beta = \limsup\limits_{n\to\infty} a_n.$

\vskip5pt

For a module $M$, put $t_0(M) = t(M)= |\top M|,$ and  $t_n(M) = t(\Omega^nM), \ (n\ge 1)$.
The numbers $t_n(M)$ are called the {\it Betti numbers} of $M$. The upper limits
$$\gamma(M) = \limsup\limits_{n\to\infty}  \root n \of {t_n(M)}$$
is called {\it the $\Omega$-growth} of $M$.
If $M$ is of finite projective dimension, then $\gamma(M) = 0$;
otherwise $\gamma(M) \ge 1.$
Note that
$$
 \gamma(M) = \limsup\limits_{n\to\infty} \root n \of {t_n(M)} = \limsup\limits_{n\to\infty} \root n \of {|\Omega^nM|}.
$$

Put  $\gamma_A = \gamma(S).$

\vskip5pt

As C. M. Ringel observed, for short local algebras, the $\Omega$-growth is a decisive invariant.

\begin{thm} \label{growth} {\rm ([RZ7, Theorem 1])} \  Let $A$ be a local algebra. Then $$\gamma(M) = \gamma(\Omega M) \le \gamma_A \le |{}_AJ|$$
for any module $M$. If
$S$ is a direct summand of $\Omega^n M$ for some $n \ge 0$,
then $\gamma(M) = \gamma_A.$ \end{thm}

\subsection{Koszul modules}
A module $M$ is {\it a Koszul module} (see J. Herzog and S. Iyengar [HI]; also [AIS]) provided a minimal projective resolution
$$\cdots \longrightarrow  P_{n+1}  \longrightarrow  P_n \longrightarrow  P_{n-1}
\longrightarrow \cdots  \longrightarrow P_0 \longrightarrow M \longrightarrow 0$$
induces for any $n\ge 0$ an exact sequence
$$0 \longrightarrow P_n/JP_n \longrightarrow JP_{n-1}/J^2P_{n-1} \longrightarrow \cdots \longrightarrow
    J^nP_0/J^{n+1}P_0 \longrightarrow J^nM/J^{n+1}M \longrightarrow 0.$$

Projective modules are always Koszul modules. If $M$ is a Koszul module over a short local algebra,
then also $\Omega M$ is a Koszul module and has Loewy length at most $2$. For a general ring, there are different definitions for a Koszul module, but
for  modules of Loewy length
at most $2$ over short local algebras, all these definitions coincide.

\begin{prop} \label{equivkoszul} {\rm ([RZ7, Proposition 4.1])} \   Let $A$ be a short local algebra and $M$ a module
of Loewy length at most $2$. The following conditions are equivalent$:$

\vskip5pt

${\rm(1)}$ \  $M$ is a Koszul module.

\vskip5pt

${\rm(2)}$ \  For every $n\ge 1$, the exact sequence
\ $0 \longrightarrow \Omega^{n}M \longrightarrow P(\Omega^{n-1}M) \longrightarrow \Omega^{n-1}M \longrightarrow 0$ induces
  an exact sequence $$0 \longrightarrow \Omega^{n}M/J\Omega^{n}M \longrightarrow
  JP(\Omega^{n-1}M)/J^2P(\Omega^{n-1}M) \longrightarrow J\Omega^{n-1}M \longrightarrow 0.$$

\vskip5pt

${\rm(3)}$ \  $\bdim \Omega^nM = (\omega^e_a)^n\bdim M$ for all $n\ge 0,$ where $\omega^e_a = \left(\begin{smallmatrix} e & -1 \cr a & 0 \end{smallmatrix}\right).$
\end{prop}

\vskip5pt

A local algebra is a {\it left Koszul algebra} if  the simple module $S$
is a Koszul module.

\begin{thm} \label{Koszul1} {\rm ([RZ7, Theorem 2])} \  Let $A$ be a short local algebra of Hilbert type $(e,a).$
If there exists a non-zero Koszul module $M$
of Loewy length at most $2$, then $A$ is left Koszul,  $a\le \frac14 e^2$ and $\gamma_A = \frac12(e+\sqrt{e^2-4a}).$
In addition, either
$\gamma(M) = \gamma_A$ or else $a > 0$ and
$\gamma(M) = \frac12(e-\sqrt{e^2-4a}\,).$\end{thm}

Theorem \ref{Koszul1} provides a generalization of the  key lemma in Lescot [Les, 3.6]:
for a commutative short local algebra $A$ with $\soc A = J^2$,
the existence of a non-zero Koszul module of Loewy length at most 2 implies that
$S$ is a Koszul module.

\begin{thm} \label{Koszul2} {\rm ([RZ7, Theorem 3])} \  Let $A$ be a short local algebra of Hilbert type $(e,a).$
Let $M$ be a non-zero module of Loewy length at most $2$ with $\gamma(M) < \gamma_A$.
Then $M$ is a Koszul module, and $\gamma(M)$ and $\gamma_A$ are positive integers with
$$
 e = \gamma(M)+\gamma_A \quad\text{and}\quad a = \gamma(M)\cdot\gamma_A.
$$
In particular, one has $0 < \gamma(M) < \frac12e < \gamma_A < e$ and
$e^2-4a = (\gamma_A-\gamma(M))^2$.
Also, $\bdim M$ is a multiple of $(1,\gamma_A)$
and $\bdim \Omega^n M = \gamma(M)^n\bdim M$ for all $n\in\Bbb N$.\end{thm}

Let $A$ be a local algebra.   An ideal $U$ of $A$ is {\it a left Conca ideal} if $U^2 = 0$ and $J^2 \subseteq JU$.
If $A$ has a left Conca ideal $U$, then $A$ is short; and since $J^2 \subseteq U$, the modules annihilated by $U$
have Loewy length at most $2$.

\begin{thm} \label{Koszul3} {\rm ([RZ7, Theorem 4])} \  Let $A$ be a short local algebra. If $A$ has a left Conca ideal $U$,
then any module annihilated by $U$ is a Koszul module$;$ in particular,  $A$ is a left Koszul algebra.\end{thm}

Theorem \ref{Koszul3} generalizes part of Theorem 1.1 of Avramov, Iyengar, \c Sega [AIS].

\subsection{Numerical descriptions of left Koszul algebras}

\begin{thm} \label{Koszul4} {\rm ([RZ7, Theorem 5, Theorem 6])} \  Given a pair $e,a$ of natural numbers, then the following are equivalent.
\vskip5pt

${\rm(1)}$ \   \ There is a short local algebra of Hilbert type $(e,a)$ which is left Koszul.

\vskip5pt

${\rm(2)}$ \    \ There is a commutative short local algebra of Hilbert type $(e,a)$ which is
left Koszul.

\vskip5pt

${\rm(3)}$ \  One has $a \le \frac14 e^2$.

\vskip5pt

Moreover, if $e = c+d, \ a = cd$ for some positive integers $c, d$,
then there are short local algebras of Hilbert type $(e,a)$ $($even commutative ones$)$
with a Koszul module with dimension vector $(1,c).$\end{thm}

\begin{thm} \label{growth2} {\rm ([RZ7, Theorem 7])} \ Let $\Cal A(a,e)$ be the class of all short local algebras of Hilbert
type $(e,a)$. If $a \le \frac14 e^2,$ then the subset
\ $\{\gamma_A\mid A\in \Cal A(e,a)\}$  of $\Bbb R$
has the least element  $\frac12(e+\sqrt{e^2-4a}).$ \end{thm}
	
\section{\bf Simple reflexive modules}

If $A$ is self-injective, then all modules are reflexive.
R. Marczinzik \cite{M3} asked whether $A$ has to be self-injective in case all the simple modules are reflexive,
and he provided several classes of algebras to support the positive answer.
C. M. Ringel [R2] answered this question in the negative, by giving an 8-dimensional wild algebra.
His example is as follows.

\begin{exm} \label{simplereflnotselfinj} \ {\rm([R2, 2.1])} \ Let $A$ be the algebra given by quiver
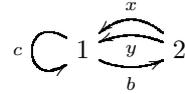
\[
\xymatrix{
		1 \ar@(ul,dl)[]_{c} \ar@/_/[r]_{b} & 2 \ar@/_1pc/[l]_{x} \ar@/_/[l]^{y}
	}
\]	
with relations \ $ cx, by, c^2, bc, xb, yb.$
Then  $A$ is not self-injective.
There are  $6$ indecomposable torsionless modules, namely the projective modules $P(1)$, $P(2)$,
the simple modules $S(1)$, $S(2)$ and the modules $P(1)/S(1)$ and $P(1)/S(2)$. Using the technique of $\mho$-quiver,
he can prove that $S(1)$ and $S(2)$ are the only indecomposable non-projective modules which are reflexive.
\end{exm}

Furthermore, [R2, 2.2] provides a generalization which is an algebra $A$ of dimension $2n+4$ with $n$ simple modules, all reflexive, but $A$ is not self-injective.
Moreover, [R2, Section 4] gives several general results on simple reflexive modules.

\vskip 5 pt

An algebra $A$ is left QF-2 if any indecomposable projective module $P$ is uniform;
and $A$ is QF-3 if the injective envelope $I(_AA)$ is projective.
Marczinzik [M3] has shown that for artin algebra $A$ satisfying one of conditions (a), (b), (c), (d) below, if
all the left simple $A$-modules are reflexive, then $A$ is self-injective. This is generalized by C. M. Ringel [R2] as follows.

\begin{thm} \label{torlselfinj} \ {\rm([R2, 5.2])} \   Let $A$ be an artin algebra satisfying one of the following conditions$:$

\vskip 5pt

{\rm(a)} \ ${\rm inj.dim.}{_{A}A} = {\rm inj.dim.}{A_A}$ and \ ${\rm proj.dim.}I(_{A}A)<\infty$ $($for example, $A$ is Gorenstein$)$.

\vskip 5pt

{\rm(b)} \ ${\rm proj.dim.}I(_{A}A) \leq 1$ $($for example, $A$ is ${\rm QF\mbox{-}3}$$)$.

\vskip 5pt

{\rm(c)} \ $A$ is left ${\rm QF}\mbox{-}2$.

\vskip 5pt

{\rm(d)}  \ Any simple right module is reflexive.

\vskip 5pt

\noindent Suppose that all the simple left $A$-modules are torsionless. Then $A$ is self-injective.
\end{thm}

More algebras $A$ are known to have the property that if all the left simple $A$-modules are reflexive then $A$ is self-injective.
For examples, local algebras ([M3, 3.14]), and algebras with radical square zero ([R2, 5.4]).

\vskip10pt

Nan Gao, \ \ \ \ \ \ \ Department of Mathematics, \ Shanghai University,  \ Shanghai 200444, \ China

\vskip5pt

Xue-Song Lu, \ \ School of Mathematical Sciences, \ Shanghai Jiao Tong University,

\hskip69pt Shanghai 200240, \ China

\vskip5pt

Pu Zhang, \ \ \ \ \ \ School of Mathematical Sciences, \ Shanghai Jiao Tong University,

\hskip69pt Shanghai 200240, \ China


\vskip20pt

\begin{thebibliography}{99}

\bibitem  [AHS]{AHS} J. Asadollahi, R. Hafezi, S. Sadeghi, On the monomorphism category of $n$-cluster tilting subcategories, Sci. China Math. 65(7)(2022), 1343-1362.

\bibitem  [A]{A} M. Auslander, Anneaux de Gorenstein et torsion en alg\'ebre commutative, S\`eminaire d'alg\'ebre
commutative dirig\'e par P. Samuel (1966-1967), tome 1, notes by: M. Mangeney, C. Peskine, L.
Szpiro, Secr\'etariat Math\'ematique, Paris, 2--69(1967).

\bibitem [AB]{AB}  M. Auslander, M. Bridger, Stable module
theory, Mem. Amer. Math. Soc. 94., Amer. Math. Soc., Providence,
R.I., 1969.

\bibitem[AR1]{AR1}   M. Auslander, I. Reiten, On a generalization of the Nakayama
  conjecture,  Proc. Amer. Math. Soc. 52 (1975), 69 - 74.

\bibitem[AR2]{AR2}  M. Auslander, I. Reiten, Applications of
contravariantly finite subcategories, Adv. Math. 86(1991), 111-152.

\bibitem[AR3]{AR3} M. Auslander, I. Reiten, Cohen-Macaulay
and Gorenstein artin algebras, In: Representation theory of finite
groups and finite-dimensional algebras, Progress in Math. vol. 95, 221-245, Birkh\"auser, Basel,
1991.

\bibitem[AR4]{AR4}   M. Auslander, I. Reiten, Syzygy modules for Noetherian rings, J. Algebra 183
(1996), 167-185.

\bibitem[ARS]{ARS} M. Auslander, I. Reiten, S. O.
Smal${\o}$, Representation Theory of Artin Algebras, Cambridge
Studies in Adv. Math. 36., Cambridge Univ. Press, 1995.

\bibitem[AGP]{AGP} L. L. Avramov, V. N. Gasharov, I. V. Peeva, Complete intersection dimension,
  Inst. Hautes Etudes Sci. Publ. Math. 86 (1997), 67 - 114.

\bibitem[AIS]{AIS} L. L. Avramov, S. B. Iyengar, L. M. \c Sega, Free resolutions
  over short local rings, J. London Math. Soc. 78(2008), 459-476.

\bibitem [AM]{AM}  L. L. Avramov, A. Martsinkovsky, Absolute, relative, and Tate cohomology of modules of
finite Gorenstein dimension, Proc. London Math. Soc. 85(3)(2002), 393-440.

\bibitem[BBD]{BBD} A. A. Beilinson, J. Bernstein, P. Deligne, Faisceaux pervers,  Ast\'erisgue 100, Soc. Math. France, Paris, 1982.

\bibitem[Bel1]{Bel1}  A. Beligiannis, The homological theory of contravariantly finite subcategories: auslander-buchweitz contexts, gorenstein categories and (co-)stabilization,
Comm. Algebra 28(10)(2000), 4547-4596.

\bibitem[Bel2]{Bel2}  A. Beligiannis, Cohen-Macaulay modules,
(co)torsion pairs and virtually Gorenstein algebras, J. Algebra
288(1)(2005), 137 - 211.

\bibitem[Bel3]{[Bel3]} A. Beligiannis, On algebras of finite Cohen-Macaulay type, Adv. Math. 226(2)(2011), 1973-2019.


\bibitem[BM]{BM} D. Bennis, N. Mahdou, Strongly Gorenstein projective, injective, and
 flat modules, J. Pure Appl. Algebra 210(2007), 437-445.

\bibitem[Bir]{Bir} G. Birkhoff,
Subgroups of abelian groups, Proc. Lond. Math. Soc. II, Ser.
38(1934), 385-401.

\bibitem [Buch]{Buch}   R.-O. Buchweitz, Maximal Cohen-Macaulay modules
and Tate cohomology over Gorenstein rings,
Unpublished manuscript, Hamburg (1987).

\bibitem[BGS]{BGS} R. O. Buchweitz, G. M. Greuel, F. O. Schreyer, Cohen-Macaulay modules on
hypersurface singularities II, Invent. Math. 88(1)(1987), 165-182.

\bibitem[Chen1]{Chen1} X. W. Chen, The stable monomorphism
category of a Frobenius category, Math. Res. Lett. 18(1)(2011),
125-137.

\bibitem[Chen2]{Chen2} X. W. Chen, Three results on Frobenius categories, Math. Z. 270(2012), no. 1-2, 43-58.

\bibitem[Chen3]{Chen3}  X. W. Chen, Algebras with radical square zero are either self-injective or CM-free,
Proc. Amer. Math. Soc. 140(1)(2012), 93-98.

\bibitem[Chr]{Chr}  L. W. Christensen, Gorenstein Dimensions, Lecture Notes in Math.
1747, Springer-Verlag, 2000.

\bibitem [CH]{CH} L. W. Christensen, H. Holm, Algebras that satisfy Auslander's condition on
  vanishing of cohomology, Math. Z. 265(2010), 21 - 40.

\bibitem[CFH]{CFH} L. W. Christensen, A. Frankild, H. Holm, On
Gorenstein projective, injective and flat dimensions-a functorial
description with applications, J. Algebra 302(1)(2006), 231-279.

\bibitem[CPST]{CPST} L. W. Christensen, G. Piepmeyer, J. Striuli,
R. Takahashi, Finite Gorenstein representation type implies simple
singularity,  Adv. Math. 218(2008), 1012-1026.

\bibitem[CV]{CV} L. W. Christensen, O. Veliche,  Acyclicity over local rings with
  radical cube zero, Illinois J. Math. 51(2007), 1439 - 1454.

\bibitem[DR]{DR}  V. Dlab, C. M.
Ringel, The module theoretical approach to quasi-hereditary
algebras, In: Representations of algebras and related topics (Kyoto,
1990), 200-224, LMS LNS 168, Cambridge
Univ. Press, 1992.


\bibitem[EJ1]{EJ1}  E. E. Enochs, O. M. G. Jenda, Gorenstein
injective and projective modules, Math. Z. 220(4)(1995), 611-633.

\bibitem[EJ2]{EJ2}  E. E. Enochs, O. M. G. Jenda, Relative
homological algebra, De Gruyter Exp. Math. 30. Walter De Gruyter
Co., 2000.

\bibitem[GKKP] {GKKP} N. Gao, J. K\"ulshammer, S. Kvamme, C. Psaroudakis,
A functorial approach to monomorphism categories II: Indecomposables, Proc. Lond. Math. Soc.(3)129(4)(2024), e12640.



\bibitem[GLS] {GLS} C. Geiss, B. Leclerc, J. Schr\"oer, Quivers with relations for symmetrizable Cartan matrices I: Foundations,
Invent. Math. 209(1)(2017), 61-158.


\bibitem[GP] {GP} I. M. Gelfand, V.A. Ponomarev, Model algebras and representations of graphs,  Funkc. anal. i. prilo$\check{z}$. 13(1979), 1-12.

\bibitem[G] {G} W. H. Gustafson, Global dimension in serial rings, J. Algebra 97 (1985), 14-16.

\bibitem[Hap1]{Hap1} D. Happel, Triangulated
categories in representation theory of finite dimensional algebras,
London Math. Soc. Lecture Notes Ser. 119, Cambridge Uni. Press,
1988.

\bibitem[Hap2]{Hap2} D. Happel, On Gorenstein algebras, in:
Representation theory of finite groups and finite-dimensional
algebras, Prog. Math. 95, 389-404, Birkha\"user, Basel, 1991.

\bibitem[HI]{HI} J. Herzog, S. Iyengar,  Koszul modules,  J. Pure Appl. Algebra  201(2005), 154 - 188.

  \bibitem[Hol]{Hol} \ H. Holm, Gorenstein homological dimensions, J. Pure Appl. Algebra
189(1-3)(2004), 167-193.

\bibitem[Hos]{Hos} M. Hoshino,  Modules without self-extensions and Nakayama's conjecture,
  Arch. Math. 43 (1982), 111 - 137.

\bibitem [HH]{HH} C. Huang, Z. Y. Huang, Torsionfree dimension of modules and self-injective
 dimension of rings, Osaka J. Math. 49(2012), 21-35.

\bibitem[HJS]{HJS} M. T. Hughes, D. A. Jorgensen, L. M. \c Sega,
  Acyclic complexes of finitely generated free modules over local ring,
  Math. Scand. 105 (2009), 85--98.

\bibitem[HSV]{HSV}
  C. Huneke, L. M. \c Sega, A. N. Vraciu,  Vanishing of Ext and
  Tor over some Cohen-Macaulay local rings, Illinois J. Math. 48(1)(2004), 295 - 317.

\bibitem [IY]{IY} O. Iyama, Y. Yoshino, Mutation in triangulated categories and rigid Cohen-Macaulay modules, Invent. Math. 172(2008), 117-168.

\bibitem[JS]{JS}  D. A. Jorgensen, L. M. \c{S}ega, Independence of the
total reflexivity conditions for modules,  Algebras and
Representation Theory 9(2)(2006), 217-226.

\bibitem[Kel1]{Kel1} B. Keller,  Chain complexes and stable categories, Manuscripta Math. 67(1990),  379 -417.

\bibitem[Kel2]{Kel2} B. Keller,  On triangulated orbit categories,  Doc. Math. 10 (2005), 551 - 581.

\bibitem[Krau]{Krau} H. Krause, Localization theory for triangulated categories, In: Triangulated categories, 161-235. London Math. Soc.
Lecture Note Ser. 375. Cambridge Univ. Press, 2010.

\bibitem[KLM1]{KLM1} D. Kussin, H. Lenzing, H. Meltzer, Nilpotent operators and weighted
projective lines, J. reine angew. Math. 685(6)(2010), 33-71.


\bibitem[KLM2]{KLM2} D. Kussin, H. Lenzing, H. Meltzer,
Triangle singularities, ADE-chains, and weighted projective lines,
Adv. Math. 237(2013), 194-251.

\bibitem[Les]{Les} J. Lescot, Asymptotic properties of Betti numbers of modules over
  certain rings, J. Pure Appl. Algebra 38 (1985), 287--298.

\bibitem[LiZ]{LiZ} Z. W. Li, P. Zhang,  A construction of Gorenstein-projective modules, J. Algebra 323(2010), 1802 - 1812.

\bibitem[LuoZ]{LuoZ} X. H. Luo, P. Zhang, Monic representations and Gorenstein-projective modules, Pacific J. Math. 264(1)(2013), 163-194.

\bibitem[L]{L} G. Lusztig, Quivers, perverse sheaves, and quantized enveloping algebras, J. Amer. Math. Soc. 4(2)(1991), 365-421.

\bibitem[M1]{M1}  R. Marczinzik, Gendo-symmetric algebras, dominant dimensions and Gorenstein
    homological algebra, arXiv:1608.04212.

\bibitem[M2]{M2}  R. Marczinzik, On stable modules that are not Gorenstein projective, arXiv:1709.01132v3.

\bibitem[M3]{M3}  R. Marczinzik, Simple reflexive modules over Artin algebras, J. Algebra Appl. 18(10)(2019), 1950193.

\bibitem[O]{O} D. Orlov, Triangulated categories of singularities and D-branes in Landau-Ginzburg
models, Proc. Steklov Inst. Math. 246(3)(2004), 227-248.




\bibitem[R1]{R1} C. M. Ringel, The Gorenstein projective modules for the Nakayama algebras, I. J. Algebra 385(2013), 241 - 261.

\bibitem[R2]{R2} C. M. Ringel,  Simple reflexive modules over finite-dimensional algebras, J. Algebra Appl. 20(9)(2021), 2150166.

\bibitem[R3]{R3} C. M. Ringel, The short local algebras of dimension $6$ with non-projective reflexive modules, Comm. Math. Stat. 11(2)(2023), 195 - 227.

\bibitem[R4]{R4} C. M. Ringel, The preprojective algebra of a quiver, In: Algebras and Modules II (Geiranger, 1996), CMS Conf. Proc. 24, 467 - 480, Amer. Math. Soc. 1998.

\bibitem[R5]{R5} C. M. Ringel, The elementary 3-Kronecker modules,  arXiv:1612.09141.

\bibitem[RS1]{RS1} C. M. Ringel, M. Schmidmeier, Submodule categories of wild representation type, J. Pure Appl. Algebra
205(2)(2006), 412-422.

\bibitem[RS2]{RS2} C. M. Ringel, M. Schmidmeier, The Auslander-Reiten translation in submodule categories, Trans. Amer. Math. Soc. 360(2)(2008), 691-716.

\bibitem[RS3]{RS3} C. M. Ringel, M. Schmidmeier, Invariant subspaces of nilpotent
operators I, J. rein angew. Math. 614(2008), 1-52.

\bibitem[RZ1]{RZ1} C. M. Ringel, P. Zhang,  From submodule categories to preprojective algebras, Math. Z. 278(1-2)(2014), 55 - 73.

\bibitem[RZ2]{RZ2} C. M. Ringel, P. Zhang,  Objective triangle functors, Sci. China Math. 58(2)(2015), 221 - 232.

\bibitem[RZ3]{RZ3} C. M. Ringel, P. Zhang, Representations of quivers over the algebra
of dual numbers, J. Algebra 475(2017), 327 - 360. Special Issue in
Memory of Prof. J. A. Green,  Edited by B. Srinivasan, M. Collins
and G. Lehrer.

\bibitem[RZ4]{RZ4} C. M. Ringel, P. Zhang,  Gorenstein-projective and
  semi-Gorenstein-projective modules, Algebra \& Number Theory 14-1(2020), 1 - 36.

\bibitem[RZ5]{RZ5}  C. M. Ringel, P. Zhang, Gorenstein-projective and
  semi-Gorenstein-projective modules II, J. Pure Appl. Algebra 224 (2020), 106248.

\bibitem[RZ6]{RZ6}  C. M. Ringel, P. Zhang, On modules $M$ such that both $M$ and $M^*$ are semi-Gorenstein-projective,
Algebr. Represent. Theory 24(2021), 1125 - 1140.

\bibitem[RZ7]{RZ7}  C. M. Ringel, P. Zhang,  Koszul modules (and the $\Omega$-growth of modules)
over short local algebras, J. Pure Appl. Algebra 225(2021), 1067772.

\bibitem[RZ8] {RZ8} C. M. Ringel, P. Zhang,  Gorenstein-projective modules over short local algebras, J. Lond.  Math. Soc. (2)106(2022), 528 - 589.

\bibitem[SWSW]{SWSW} S. Sather-Wagstaff, T.Sharif, D. White, Stability of Gorenstein categories, J. Lond. Math. Soc. 77(2)(2008) 481-502.

\bibitem[S1]{S1} D. Simson, Representation types of the
category of subprojective representations of a finite poset over
$K[t]/(t^m)$ and a solution of a Birkhoff type problem, J. Algebra
311(2007), 1-30.

\bibitem[S2]{S2} D. Simson, Representation-finite Birkhoff type problems for
nilpotent linear operator,  J. Pure Appl. Algebra 222(2018), 2181-2198.

\bibitem[V]{V} M. Van den Bergh, Non-commutative crepant resolutions, In: The legacy of Niels Henrik Abel, Berlin: Springer-Verlag, 2004, 749-770.

\bibitem[Y]{Y}  Y. Yoshino, A functorial approach to modules of $G$-dimension zero, Illinois J. Math. 49(3)(2005), 345-367.

\bibitem[Z1]{Z1} P. Zhang, Monomorphism categories, cotilting theory, and Gorenstein-projective modules, J. Algebra
339(2011), 180-202.

\bibitem[Z2]{Z2} P. Zhang, Monic modules and semi-Gorenstein-projective modules, J. Pure Appl. Algebra 227(2023), 107181.

\bibitem[ZX]{ZX}  P. Zhang, B. L. Xiong,  Separated monic representations II: Frobenius subcategories and RSS equivalences,  Trans. Amer. Math. Soc. 372(2)(2019), 981-1021.

\end{thebibliography}
\end{document}